\documentclass[journal]{IEEEtran}

\usepackage{mystyle}

\begin{document}

\title{Plugging Weight-tying Nonnegative Neural Network into Proximal Splitting Method: Architecture for Guaranteeing Convergence to Optimal Point}
\author{Haruya~Shimizu,
        Masahiro~Yukawa,~\IEEEmembership{Senior Member,~IEEE,}
\thanks{
This work was supported in part by JSPS Grants-in-Aid (23K22762).
%
}
\thanks{
The authors are with the Department of Electronics and Electrical Engineering, Keio University, Yokohama, Kanagawa 223-8522, Japan (e-mail: yukawa@elec.keio.ac.jp).
}}


\maketitle


\begin{abstract}
We propose a novel multi-layer neural network architecture
that gives a promising neural network empowered optimization approach
to the image restoration problem.
The proposed architecture is motivated by the recent study of monotone
Lipschitz-gradient (MoL-Grad) denoiser (Yukawa and Yamada, 2025) which
establishes an ``explainable'' plug-and-play (PnP) framework
in the sense of disclosing the objective minimized.
The architecture is derived from the gradient of a superposition of
functions associated with each layer, having the weights in the encoder
and decoder tied with each other.
Convexity of the potential, and thus monotonicity of its gradient
(denoiser),
is ensured by restricting ourselves to nonnegative weights.
Unlike the previous PnP approaches with theoretical guarantees,
the denoiser is free from constraints on the Lipschitz constant
of the denoiser.
Our PnP algorithm employing the weight-tying nonnegative neural network
converges to a minimizer of the objective involving an ``implicit''
weakly convex regularizer induced by the denoiser.
The convergence analysis relies on an efficient technique to preserve
the overall convexity even in the ill-conditioned case where
the loss function is not strongly convex.
The simulation study shows the advantages of the
 Lipschitz-constraint-free nature of the proposed denoiser
in training time as well as deblurring performance.
\end{abstract}
\begin{IEEEkeywords}
Plug and Play (PnP), proximal operator, weakly convex function, image restoration, convergence guarantee
\end{IEEEkeywords}

\IEEEpeerreviewmaketitle

\section{Introduction}
\IEEEPARstart{T}{his} article presents
a powerful mathematically-rigorous approach
to the image restoration problem
exploiting the proximal splitting method
assisted by neural network denoisers.
Our approach relies on the recently established framework of
the monotone Lipschitz-gradient (MoL-Grad) denoiser
\cite{yukawa2025_molgrad} which characterizes
the class of single-valued proximity (s-prox) operators of weakly convex
functions.
%
This work goes one step further towards
constructions of high performance MoL-Grad denoisers
by studying the architecture of multi-layer neural network.

The image restoration task is typically cast as the following minimization problem:
\begin{equation}
\label{inverse_problem_of_IR}
\minimize_{\bm{x}}\hspace{5pt} \frac{1}{2}\|\bm{A}\bm{x}-\bm{y}\|^2+\lambda g(\bm{x}),
\end{equation}
where $\bm{y}=\bm{A}\bm{x}^*+\bm{\varepsilon}$ is the degraded version
of the original image $\bm{x}^*$ to be restored
with the knowledge of the degradation matrix
$\bm{A}$
in the presence of the Gaussian noise vector
$\bm{\varepsilon}$,
and
$g$ is the regularizer accommodating prior knowledge about the images.
As the classical convex regularizers such as the $\ell_1$ norm,
the total variation (TV) norm \cite{RUDIN1992_TV},
and the weighted norm involving the wavelet basis
\cite{KHORAMIAN2012_penaltynorm},
tend to cause nonnegligible estimation biases,
many non-convex regularizers have been proposed in the context of image
processing tasks
\cite{jian_nonconvexTV, Rowl,Ochs_nonconv_regu, Ben_nonconvexTV}.
Instead of using regularizers explicitly,
the plug-and-play (PnP) methods \cite{ven2013_PnP}
leverage prior information encoded in the neural network denoiser which
replaces the proximity operator in the splitting algorithms.
The original work of PnP \cite{ven2013_PnP} supposed to employ
the classical denoisers such as K-SVD \cite{K_SVD_origin, K_SVD},
block-matching and 3D filtering (BM3D) \cite{BM3D}, and TV
\cite{TV_for_Pnp}.
Further developments have explored the use of learned neural network denoisers,
such as
denoising convolutional neural network (DnCNN) \cite{DnCNN},
IRCNN \cite{IRCNN},
fast and flexible denoising convolutional neural network (FFDNet)
\cite{FFDNet, FFD_for_PnP}, and
DRUNet (combining ResNet and U-Net) \cite{Zang2022_drunet},
as well as approaches based on diffusion models \cite{diffusion_PnP_IR}.
Despite the remarkable enhancement of restoration performance,
PnP methods lack ``explainability'' in general because
the regularization is induced by the denoiser
(not directly to the function $g$)
unlike the typical optimization approaches.
The explainability aspect remains an active research area of
paramount importance.

%
%
Most of the previous studies in this direction assume
the nonexpansiveness of the denoisers
\cite{sreehari2016, romano2016_red,
sun2019, teodoro2019, pesquet2021_MMO, Nair2021}.
For instance, \cite{pesquet2021_MMO} and \cite{sun2021} have studied
the firmly nonexpansive denoiser and have shown that
the PnP algorithm converges to
a zero of an operator induced by the denoiser.
In \cite{Nair2021}, the averaged nonexpansive denoiser has been studied
under other assumptions on the denoiser, and
convergence analysis of the PnP algorithm has been provided.
Regularization by denoising (RED) \cite{romano2016_red}
relies on the assumptions of
nonexpansiveness, local homogeneity, and Jacobian symmetry,
and its extension called RED-PRO \cite{RED_pro}
reformulates it using the projection onto the fixed-point set of
the denoiser, serving as a bridge between PnP and RED.
Motivated by the fact that the conditions of RED
are unsatisfied by deep denoisers \cite{reehost2019},
the gradient-step (GS) denoiser has been proposed
\cite{gradient_step_origin, hurault2022_gradient}.
In \cite{hurault2022proximal}, it has shown
that the GS denoiser is actually the proximity operator of a certain regularizer
if the residual $\Id-D$ of denoiser $D \ (=\nabla h)$
is contraction
(i.e., Lipschitz constant strictly smaller than unity),
where the function $h$ is a potential.
Unfortunately, however, typical CNN denoisers
are non-nonexpansive (and fail to satisfy even the relaxed variants
\cite{RED_pro, hurault2022proximal}).
As such, CNN-based approaches typically incorporate additional
penalty terms into the training objective to enforce
the condition on the Lipschitz constant
\cite{ryu2019_nonexptrain_by-norm, pesquet2021_MMO, hurault2022proximal}.
Imposing such Lipschitz constraints on CNNs, however,
often causes performance degradation as well as
incurring a substantial computational burden
due to the repetitive use of power iterations for training the network.\footnote{
While the Lipschitz constant constraints have been considered for technical reasons to guarantee convergence of algorithms, the Lipschitz constant is also related closely to robustness to adversarial perturbations \cite{tsuzuku2018, pesquet_icassp20}.
%
Indeed, there is an intrinsic tradeoff between accuracy and robustness in general, and the Lipschitz constant tends to be large for attaining high accuracy \cite{pesquet_icassp20}.
}

The MoL-Grad denoiser \cite{yukawa2025_molgrad} has emerged
as a powerful mathematical framework to address this issue.
Unlike the aforementioned approaches,
MoL-Grad has no Lipschitz constraint, and therefore
it is more expressive than those with the nonexpansiveness constraint
(or its relaxed version).
%
%
The backbone of the theory of MoL-Grad is the fact that
the s-prox operator of weakly convex
functions can be characterized by monotone Lipschitz-continuous gradient
operators,
or, equivalently, by gradient operators $\nabla \psi$ of smooth convex
functions $\psi$ (see Fact \ref{fact:def_molgrad}).
By virtue of an extended version of Moreau's decomposition
\oldcite[Proposition 2]{yukawa2025_molgrad},
it has been shown that,
adopting the MoL-Grad denoiser $\nabla \psi$ in the PnP method,
the primal-dual splitting type algorithm, as well as the forward-backward splitting type algorithm,
converges to a minimizer of a certain cost function that involves
the ``implicit'' regularizer $\varphi := \psi^* - (1/2) \|\cdot\|^2$.
Although a simple weight-tying two-layer neural network has been shown to be
a MoL-Grad denoiser,
exploring more practical multi-layer neural networks under the MoL-Grad framework
is an attractive open issue of wide interest
especially when one considers its application to image restoration.

In this paper, we study a multi-layer neural network architecture that yields
a MoL-Grad denoiser.
The proposed architecture has a structure of auto-encoder, and
since the denoiser needs to possess a symmetric Jacobian,
the weights in the encoder and decoder are tied with each other.
We also restrict our current attention to nonnegative weights
to ensure monotonicity of the denoiser.
Our weight-tying multi-layer nonnegative neural network is shown to be a MoL-Grad
denoiser under certain assumptions.
As the proposed neural network  denoiser is associated ``implicitly''
with a weakly convex regularizer,
 the overall cost function could be nonconvex when used to solve
ill-conditioned problems such as image restoration.
To circumvent this difficulty, we present a technique to
restrict the weak convexity of the implicit regularizer to a certain subspace
so that the overall convexity is preserved.
Thanks to this restriction technique, the proposed algorithm converges to
a global minimizer of the cost function involving
the restricted counterpart of the implicit regularizer.
The simulation results show the remarkable advantages of the proposed
denoiser in terms of the training time as well as its state-of-the-art
performance in application to image deblurring.
The remainder of the paper is organized as follows.
After the notation and definitions in Section~\ref{sec:preliminaries},
Section~\ref{sec:proposed} presents the proposed weight-tying
multi-layer nonnegative neural network architecture, which is shown to give
a MoL-Grad denoiser under some conditions.
Section~\ref{sec:algorithm} then presents the proposed algorithm and its
convergence analysis, and
Section~\ref{sec:simulation} shows the simulation results,
followed by conclusion in Section~\ref{sec:conclusion}.

%
%
\section{Preliminaries}
\label{sec:preliminaries}
We first present notations and definitions, and we then introduce the MoL-Grad denoiser.
\subsection{Notation}
Let $\mathbb{R}$, $\mathbb{R}_+$, $\mathbb{R}_-$, and $\mathbb{N}$ be the sets of real numbers, nonnegative real numbers, non-positive real numbers, and nonnegative integers, respectively.
The superscripts $(\cdot)^\mathsf{T}$ stand for transposition.
Let $\Id$ be the identity operator.
For a vector $\bm{x} \in \mathbb{R}^n$, the $x_i$ denotes the $i$th component, and the $\ell_2$ norm is defined by $\|\bm{x}\| := \left(\sum_{i=1}^{n}|x_i|^2\right)^{1/2}$.
For a matrix $\bm{A} \in \mathbb{R}^{n \times m}$, let $[\bm{A}]_i \in \mathbb{R}^{1 \times n}$ denote its $i$th row, and $[\bm{A}]_{i, j}$ denote its $(i, j)$ entry, and the spectral norm is defined by $\|\bm{A}\|:=\max_i\{\sigma_i(\bm{A})\}$ where $\sigma_i(\bm{A})$ denotes the $i$th singular value of $\bm{A}$.
A function $f:\mathbb{R} \to \mathbb{R}$ is non-decreasing if $x \leq y$ implies $f(x) \leq f(y)$.
A function $f$ is convex if $f(a\bm{x} + (1-a)\bm{y}) \leq af(\bm{x}) + (1-a)f(\bm{y})$ for every $\bm{x}, \bm{y} \in \mathbb{R}^n$ and $a \in (0,1)$.

The function $f$ is $\rho$-strongly convex if $f - (\rho/2) \|\cdot\|^2$ is convex, and $\rho$-weakly convex if $f + (\rho/2) \|\cdot\|^2$ is convex for $\rho > 0$.
The set of proper lower semicontinuous convex functions from $\mathbb{R}^n$ to $\mathbb{R}\cup\infty$ is denoted by $\Gamma_0(\mathbb{R}^n)$.
Given a function $f\in\Gamma_0(\mathbb{R}^n)$, its Fenchel conjugate is defined by $f^*: \bm{z} \mapsto \text{\normalfont sup}_{\bm{x}\in\mathbb{R}^n} \big( \langle \bm{x}, \bm{z} \rangle-f(\bm{x}) \big)$.
For $\kappa>0$, an operator $T:\mathbb{R}^n\to\mathbb{R}^n$ is $\kappa$-Lipschitz continuous if $\|T(\bm{x})-T(\bm{y})\|\leq\kappa\|\bm{x}-\bm{y}\|$ for every $\bm{x}, \bm{y} \in\mathbb{R}^n$.
In particular, if $\kappa= 1$, $T$ is a nonexpansive operator.
Moreover, $T$ is $\alpha$-averaged nonexpansive for an $\alpha\in(0, 1)$ if there exists a nonexpansive operator $Q$ such that $T=\alpha Q+(1-\alpha)\Id$.
A $1/2$-averaged nonexpansive operator $T$ is called firmly nonexpansive.
For $\eta>0$, $T$ is $\eta$-cocoercive if $\eta T$ is firmly nonexpansive.
Let the operator $T:\mathbb{R}^n \to \mathbb{R}^m$ and the function $f:\mathbb{R}^m\to\mathbb{R}$ be differentiable at every $\bm{x}\in\mathbb{R}^n$ and $\bm{y}\in\mathbb{R}^m$, respectively.
Then, the gradient vector $\nabla (f\circ T)(\bm{x})$ of the composite function $f \circ T$ at $\bm{x}$ is given by
\begin{equation}
    \nabla (f\circ T)(\bm{x})=\mathrm{J}_T(\bm{x})^\mathsf{T} \nabla f (T(\bm{x})),
\end{equation}
where $\mathrm{J}_T(\bm{x})$ is the Jacobian of $T$ at $\bm{x}$.

\begin{figure*}[t]
    \centering
    \includegraphics[width=1\linewidth]{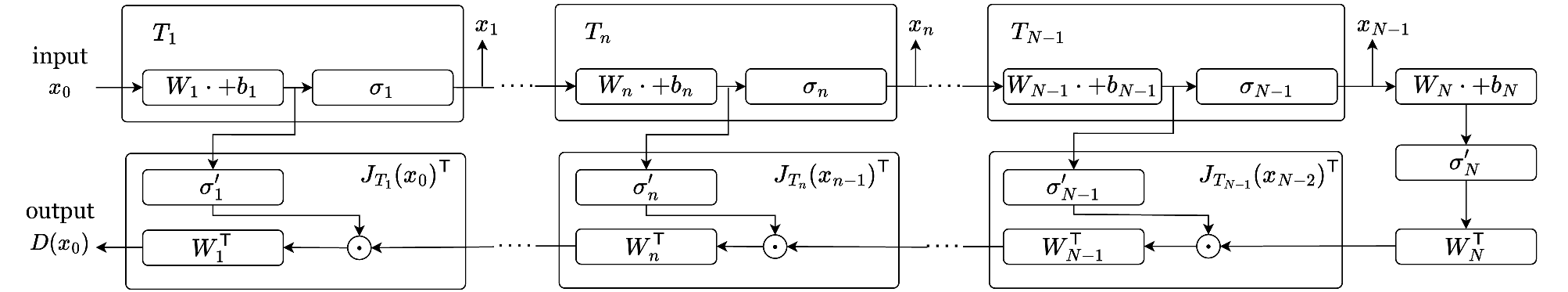}
    \caption{The architecture of the denoiser $D$.}
    \label{fig:model_Jacobian}
\end{figure*}

\subsection{MoL-Grad Denoiser}
\begin{definition}{Single-valued proximity operator}
\label{def:sprox}
    For a proper function $f:\mathbb{R}^n\to(-\infty, +\infty]$,
the s-Prox operator
$\text{\normalfont s-Prox}_{\gamma f}: \mathbb{R}^n\to\mathbb{R}^n$
of index $\gamma>0$ is defined by \cite{yukawa2025_molgrad}
\begin{align}
    \text{\normalfont s-Prox}_{\gamma f}:\bm{x}\mapsto \argmin_{\bm{y}\in\mathbb{R}^n}\left( f(\bm{y})+\frac{1}{2\gamma}\|\bm{x}-\bm{y}\|^2 \right),
\end{align}
whenever $f+(1/2\gamma)\|\bm{x}-\cdot\|^2$ has a unique minimizer for every $\bm{x}\in\mathbb{R}^n$.
\end{definition}
\begin{fact}{MoL-Grad denoiser\oldcite[Theorem 1]{yukawa2025_molgrad} for Euclidean case}
\label{fact:def_molgrad}
Let $D:\mathbb{R}^n\to\mathbb{R}^n$ be a nonlinear operator. For every $\beta\in(0,1)$, the
following two conditions are equivalent.
\begin{enumerate}[label=\normalfont(\alph*)]
    \item $D=\text{\normalfont s-Prox}_\varphi$ for some $\varphi:\mathbb{R}^n\to(-\infty, \infty]$ such that $\varphi+[(1-\beta)/2]\|\cdot\|^2 \in\Gamma_0(\mathbb{R}^n)$.
    \item $D$ is a MoL-Grad denoiser, i.e.,
it satisfies the following two conditions jointly:
    %
{\rm (i)} $D$ is $\beta^{-1}$-Lipschitz continuous, and
{\rm (ii)}  $D=\nabla\psi$ for a Fr\'{e}chet differentiable convex function $\psi$.
    Conditions {\rm (i)}  and {\rm (ii)}  are collectively referred to as condition $\clubsuit$.
    \label{fact:condition_mol}
\end{enumerate}
If condition (a), or equivalently condition (b) is satisfied,
the identity $\varphi=\psi^*-(1/2)\|\cdot\|^2$ holds.
\end{fact}
%
%
\section{Proposed Weight-tying Architecture of Nonnegative Neural
 Network Denoiser}
\label{sec:proposed}
Neural network denoisers are not MoL-Grad (i.e., unsatisfy condition
$\clubsuit$) in general.
We present a multi-layer neural network architecture
that gives a MoL-Grad denoiser under appropriate conditions.
%
%
%
After presenting the proposed architecture where the weights are tied between
encoder and decoder,
we derive a sufficient condition for the resultant weight-tied neural network
to be a MoL-Grad denoiser.

%

%
\subsection{Network Architecture of denoiser}

We first define multi-layer neural network which will then be used to
derive the network architecture.

\begin{definition}{Multiple layer neural network $\bm{T_{n:m}}$}
  \label{def_T}
  \\ {\normalfont (a)} For every $n\in\{1, 2, \cdots, N\}$, the single layer neural network $T_n:\mathbb{R}^{d_{n-1}}\to\mathbb{R}^{d_{n}}$ is defined by
  \begin{align}
    T_{n}(\bm{x}) &:=\left[ t_{n, 1}(\bm{x}), t_{n,2}(\bm{x}), \cdots, t_{n, d_n}(\bm{x})\right]^\mathsf{T}\\
    &:= \bm{\sigma}_n ( \bm{W}_{n} \bm{x}+\bm{b}_{n} ),
\end{align}
where $\bm{W}_{n} \in \mathbb{R}^{d_{n} \times d_{n-1}}$ is the weight matrix, and $\bm{b}_{n} \in \mathbb{R}^{d_{n}}$ is the bias vector, and
\begin{align}
     \bm{\sigma}_n&:\mathbb{R}^{d_n}\to\mathbb{R}^{d_n}\notag\\
                        &:\bm{x}\mapsto\left[
          \sigma_{n, 1} (x_1), \sigma_{n, 2}(x_2),
          \cdots,
          \sigma_{n, d_n} (x_{d_n}) \right]^\mathsf{T} \label{eq:def_sigma_f}
  \end{align}
is the componentwise activation operator with $\sigma_{n, i}: \mathbb{R}\to\mathbb{R}$.
\\
{\normalfont (b)} Multi-layer neural network $T_{n:m}:\mathbb{R}^{d_{m-1}}\to\mathbb{R}^{d_{n}}$ from the $m$th layer to the $n$th layer is defined by
  \begin{align}
      T_{n:m}:=T_n \circ T_{n-1} \circ \cdots \circ T_m \label{eq:def_T_n:m}
  \end{align}
  for $1\leq m \leq n \leq N$.
\end{definition}
We shall define $D$ using the $N$-layer neural network $T_{N:1}$.
\vspace{-10pt}
\begin{definition}{Denoiser $\bm{D}$}
\label{def_D}
The denoiser $D$ is defined by
\begin{equation}
    D:=\nabla\psi \quad \mbox{for} \quad\psi:= \sum_{i=1}^{d_N} t_i.
\end{equation}
%
Here, the functions $t_i:\mathbb{R}^{d_0}\to\mathbb{R}$ are defined by
\begin{align}
    T(\bm{x})&:\hspace{2pt}=\left[ t_1(\bm{x}), t_2(\bm{x}), \cdots,
 t_{d_N}(\bm{x})\right]^\mathsf{T} := T_{N:1}(\bm{x}).
\end{align}
Using the chain rule, the output of $D$ for an input vector
$\bm{x}_0 \in\mathbb{R}^{d_0}$ is given by
\begin{align}
& \hspace{-.6em}  D(\bm{x}_0) = \nabla \left( \sum^{d_N}_{i=1} t_i  \right) (\bm{x}_0)
        = \mathrm{J}_{T_{1}}(\bm{x}_{0})^\mathsf{T} \cdots \ \mathrm{J}_{T_{N-1}}(\bm{x}_{N-2})^\mathsf{T}  \nonumber\\
 & \hspace{.2em}  \circ \left[\bm{W}_{N}^\mathsf{T} \circ \bm{\sigma}_N' \circ (\bm{W}_{N} \cdot + \bm{b}_{N}) \right] \circ T_{N-1} \circ \cdots \circ T_{1} (\bm{x}_0),\label{eq:express_D}
\end{align}
%
where
\begin{align}
    \bm{x}_n:=&~T_n(\bm{x}_{n-1})=T_{n:1}(\bm{x}_0),~~~
n\in\{1, 2, \cdots, N\},\label{eq:def_x_n}\\
\hspace*{-.4em}    \bm{\sigma}_n'(\bm{x}) :=&~
\! \left[ \sigma'_{n, 1} (x_1), \sigma_{n, 2}' (x_2),\cdots, \sigma'_{n,d_n}(x_{d_n})\right]^\mathsf{T} \nonumber\\
    =&~
\! \left[ \frac{\mathrm{d} \sigma_{n, 1} (x_1)}{ \mathrm{d} x_1},
 \frac{\mathrm{d} \sigma_{n, 2}(x_2)}{\mathrm{d} x_2}, \cdots,
 \frac{\mathrm{d}\sigma_{n, d_n}(x_{d_n})}{\mathrm{d}
 x_{d_n}}\right]^\mathsf{T} \!\!\!, \label{eq:def_sigma_prime}
\end{align}
and
$\mathrm{J}_{T_{n}}(\bm{x}_{n-1}):=\mathrm{diag}(\bm{\sigma}'_n(\bm{W}_n\bm{x}_{n-1}+\bm{b}_n))\bm{W}_n$ is the Jacobian matrix of $T_{n}$ at $\bm{x}_{n-1}$.
We can then write
\begin{equation}
  \mathrm{J}_{T_{n}}(\bm{x}_{n-1})^\mathsf{T} \bm{z} = \bm{W}_{n}^\mathsf{T} (\bm{\sigma}_n'( \bm{W}_{n} \bm{x}_{n-1} + \bm{b}_n) \odot \bm{z}),~\bm{z}\in\mathbb{R}^{d_n}, \label{eq:express_jaco}
\end{equation}
where $\odot$ represents the Hadamard product between vectors.
The denoiser $D$ is the weight-tying neural network consisting of the encoder $\bm{\sigma}_N' \circ (\bm{W}_N \cdot \bm{b}_N) \circ T_{N-1} \circ \cdots \circ T_{1}$ and the decoder $\mathrm{J}_{T_{1}}(\bm{x}_{0})^\mathsf{T} \cdots \mathrm{J}_{T_{N-1}}(\bm{x}_{N-2})^\mathsf{T} \bm{W}_N^\mathsf{T}$ (see Fig.~\ref{fig:model_Jacobian}).
\end{definition}
%

%
%
\subsection{Conditions for $D$ to be MoL-Grad  denoiser}
We present a sufficient condition to guarantee convexity of $\bm{\psi}$ and Lipschitz continuity of $\nabla \psi$.
First, we introduce lemmas to derive a condition to guarantee these properties.
\begin{lemma}{}
\label{claim_conv}
 Define $T:\mathbb{R}^p\to\mathbb{R}^{q}:x\mapsto[t_1(\bm{x}), t_2(\bm{x}), \cdots, t_q(\bm{x})]^\mathsf{T}$, and assume that $t_i$ is convex for $i\in\{1, 2, \cdots, q\}$.
 In addition, let $\sigma:\mathbb{R}\to\mathbb{R}$ be a non-decreasing convex function, and $\bm{w}\in\mathbb{R}_{+}^{q}$ be a nonnegative vector.
 Then, the function $\sigma(\bm{w}^\mathsf{T}T(\cdot)+b)$ is convex for a constant $b\in\mathbb{R}$.
 \myproof
    See Appendix \ref{app:proof_pf_lemma1}. \qed
\end{lemma}
\begin{lemma}{}
\label{memo_lip2}
Let $\bm{W}\in\mathbb{R}^{p\times q}$, and let $T:\mathbb{R}^s\to\mathbb{R}^q: \bm{x}\mapsto[t_1(\bm{x}), t_2(\bm{x}), \cdots, t_q(\bm{x})]^\mathsf{T}$ and $R:\mathbb{R}^s\to\mathbb{R}^q: \bm{x}\mapsto[r_1(\bm{x}), r_2(\bm{x}), \cdots, r_q(\bm{x})]^\mathsf{T}$.
In addition, define the operator $F:\mathbb{R}^s \to \mathbb{R}^p$ and the function $f_{i}:\mathbb{R}^s\to \mathbb{R}$ as follows:
\begin{align*}
    F(\bm{x}) :=\left[ f_{1}(\bm{x}), f_{2}(\bm{x}), \cdots, f_{p}(\bm{x})\right]^\mathsf{T}
    :=  \bm{W}  \left[T(\bm{x})\odot  R(\bm{x})\right].
\end{align*}
%
Assume that the operators $T$ and $R$ are Lipschitz continuous. Furthermore, the functions $t_i$ and $r_i$ are bounded for $i\in\{1,2,\cdots, q\}$, i.e., there exist positive constants $M_t$ and $M_r$ such that $|t_i(\bm{x})|\leq M_t$ and $|r_i(\bm{x})|\leq M_r$ for every $\bm{x}\in\mathbb{R}^{s}$.
Then, the operator $F$ is Lipschitz continuous, and the function $f_i$ are bounded for $i\in\{1,2,\cdots, s\}$.
\myproof
See Appendix \ref{app:proof_pf_lemma2}. \qed
\end{lemma}
The following theorem shows that $D$ is a MoL-Grad denoiser under appropriate conditions.
\begin{theorem}{$\bm{D}$ is MoL-Grad denoiser}
\label{theorem_mol_NN}
Let $T$, $t_i$, and $D$ be defined as in Definition \ref{def_D}. \\
\noindent{\normalfont (a)} Assume the following two conditions:
  \begin{enumerate}[label=\normalfont(K.\arabic*), ref=\normalfont(K.\arabic*), leftmargin=35pt]
    \item The activation function $\sigma_{n, i}$ (see (\ref{eq:def_sigma_f})) is non-decreasing, differentiable, and convex, for every layer $n$ and component $i$.\label{assumption_conv_sigma}
    \item The weight matrix $\bm{W}_{n}$ is nonnegative for every $n \geq 2$.\label{assumption_conv_W}
  \end{enumerate}
  Then, the function $t_i$ is a differentiable convex function for every $i$.\\
\noindent{\normalfont (b)} Assume the following two conditions:
     \begin{enumerate}[label=\normalfont(L.\arabic*), ref=\normalfont(L.\arabic*), leftmargin=35pt]
    \item The activation function $\sigma_{n, i}$ is Lipschitz continuous for every $n$ and for every $i$.
    \label{assumption_lip_sigma}
    \item $\bm{\sigma}'_{n}$ (see (\ref{eq:def_sigma_prime})) is Lipschitz continuous for every $n$.
    \label{assumption_lip_sigmaprime} 
  \end{enumerate}
  Then, the denoiser $D$ is Lipschitz continuous.\\
\noindent{\normalfont (c)} Assume that \ref{assumption_conv_sigma}, \ref{assumption_conv_W}, \ref{assumption_lip_sigma} and \ref{assumption_lip_sigmaprime} are satisfied.
Then, $D$ is a MoL-Grad denoiser, i.e.,~$D=\text{s-Prox}_\varphi$ for
 some $(1-L_D^{-1})$-weakly convex function $\varphi$, under the
 assumption that $L_D > 1$ is a Lipschitz constant of $D$ (see Fact~\ref{fact:def_molgrad}).
\myproof
See Appendix \ref{app:proof_of_theorem1}. \qed
\end{theorem}

\subsection{Remarks on Theorem \ref{theorem_mol_NN}}
%

%
%
%
%
\begin{remark}{A stable way of imposing nonnegativity in network training}
\label{remark:nonneg}
Enforcing assumption \ref{assumption_conv_W} of Theorem
 \ref{theorem_mol_NN} from the beginning of training tends to cause
 instability.
To avoid this issue, the network is trained according to the following loss:
\begin{equation}
  \mathcal{L}(\bm{\theta})=
   \|D_{\bm{\theta}}(\bm{x}^*+\bm{\varepsilon})-\bm{x}^*\|^2 +
   \alpha_\mathcal{L} \sum_i  \left( \max\{0, -\theta_{w_i}\} \right)^2,
\label{eq:def_loss}
\end{equation}
where the second term encourages nonnegativity.
Here, $\bm{\theta}$ is the learnable parameters including the weights $\theta_{w_i}$, and $D_{\bm{\theta}}$ is the denoiser parameterized by $\bm{\theta}$.
The penalty in \eqref{eq:def_loss} is implemented using a quadratic barrier function \cite{Numerical_Optimization}, following the approach in \cite{nguyen_nonnegativeRBM}, \cite{hosseini_nonnegative_aut}.
In our preliminary experiments, minimizing \eqref{eq:def_loss} encourages the weights to be nonnegative softly.
Those weights which become negative-valued (if any) are enforced to be zero by post-processing.
\end{remark}

\begin{remark}{A limitation of nonexpansive denoiser in sparse modeling}
  \label{remark:nonexpansive_limitation}
Let us consider a sparse modeling problem where
the task is to find a sparse vector $\bm{x}$ such that
$\bm{Ax}$ for a given matrix $\bm{A}$ is sufficiently close to the given
output $\bm{y}$.
For this task, an ideal shrinkage operator will have the following two
properties:
(i) those components of small magnitude are attracted to zero, and
(ii) those of sufficiently large magnitudes are preserved.
To build a continuous operator having this ideal properties,
those components of moderate magnitudes need to take intermediate values
so that the Lipschitz constant exceeds one;
the reader may consider one dimensional case to see this easily.
This observation suggests a limitation of
nonexpansive shrinkage operators
(such as soft shrinkage), and thus that of nonexpansive denoisers.
Indeed, the proposed denoiser free from Lipschitz constant bounds
exhibits better performance in image restoration compared to
the nonexpansive denoisers, as will be shown in Section
 \ref{sec:simulation}.
\end{remark}

\begin{remark}{Case of $\bm{L_D\leq1}$ in Theorem~\ref{theorem_mol_NN}(c)}
When $L_D \leq 1$, $D$ is Moreau's proximity operator of convex function \cite{convex_analysis}.
Even in such a case, Theorem~\ref{theorem_mol_NN}(c) holds for another constant $L'_D > 1 \ (\geq L_D)$ because $D$ is also $L'_D$-Lipschitz continuous
by definition of Lipschitz continuity.
\end{remark}
\subsection{An Extension with Skip Connection}
The results of Theorem~\ref{theorem_mol_NN} can be extended by using skip connections, which have been widely used in neural network denoisers.
\begin{figure}
    \centering
    \includegraphics[width=1\linewidth]{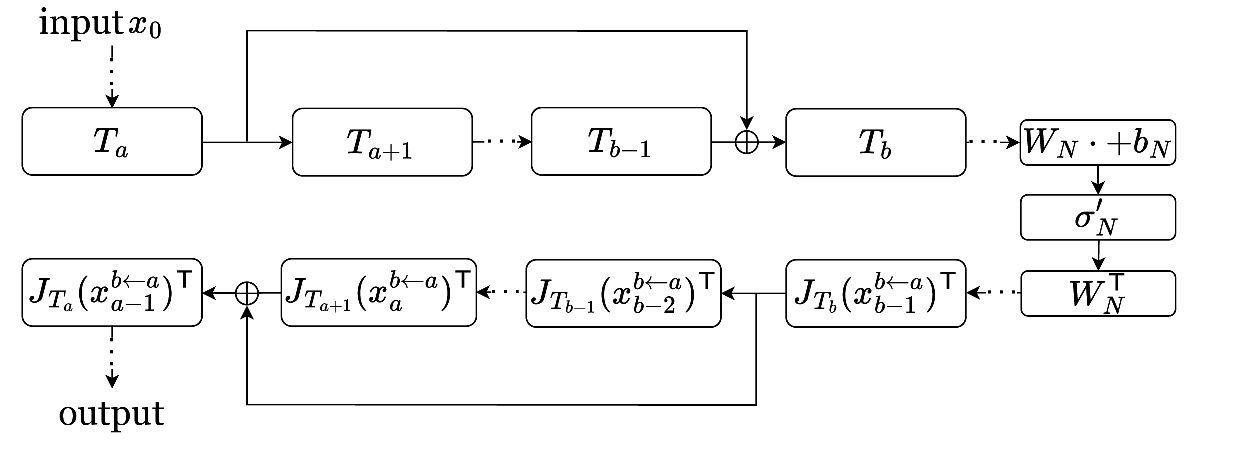}
    \caption{The architecture of $D^{b\leftarrow a}$.}
    \label{fig:structure_D_skip}
\end{figure}
\begin{theorem}{$D^{b\leftarrow a}$ with skip connection is MoL-Grad denoiser}
\label{theorem_skip}
    For $a\in\{1,2,\cdots, N-2\}$ and $b\in\{3,4,\cdots, N\}$ such that
 $b-a\geq 2$, define the denoiser involving skip connections by (see
 Fig.~\ref{fig:structure_D_skip})
    \begin{equation}
        D^{b\leftarrow a}:=\nabla \left(\sum_{i = 1}^{d_N} t_i^{b\leftarrow a}\right),
    \end{equation}
where
    \begin{align}
    T^{b\leftarrow a} (\bm{x}) :\hspace{-3pt}&=\left[ t^{b\leftarrow a}_1 (\bm{x}), t^{b\leftarrow a}_{2} (\bm{x}), \cdots, t^{b\leftarrow a}_{d_N} (\bm{x})\right]^\mathsf{T}\\
    :\hspace{-3pt}&= T_{N:b}\circ(\Id+T_{b-1:a+1})\circ T_{a:1} (\bm{x})\\
    & = T_{N:b} \circ( T_{a:1}+T_{b-1:1})(\bm{x}). \notag
    \end{align}
    Here, $T^{b\leftarrow a}$ denotes a modified operator of $T$ with
 the $a$th layer output connected to the $b$th layer input by
 skipping the intermediate layers.\footnote{
    Note that $T^{b\leftarrow a}=T_{N:b}\circ (T_{a:1}+T_{b-1:1})\neq T_{N:b}\circ T_{a:1}+T_{N:b}\circ T_{b-1:1}$ because of the nonlinearity of $T_{N:b}$.
    }
    The notation $T_{n:m}$ follows the definition in \eqref{eq:def_T_n:m}.
    Assume that assumptions \ref{assumption_conv_sigma}, \ref{assumption_conv_W}, \ref{assumption_lip_sigma} and \ref{assumption_lip_sigmaprime} of Theorem~\ref{theorem_mol_NN} are satisfied.
    Then, {\normalfont (a)\hspace{3pt}} $t^{b\leftarrow a}_{i}$ is convex for every $i$, and {\normalfont (b)\hspace{3pt}} $D^{b\leftarrow a}$ is a MoL-Grad denoiser.
\myproof
    See Appendix \ref{app:proof_of_theorem2}. \qed
\end{theorem}
\section{PnP Algorithm for Image Restoration}
\label{sec:algorithm}
We present an ``explainable'' PnP algorithm employing the MoL-Grad denoiser $D$ that is trained under the conditions of Theorem~\ref{theorem_mol_NN} or Theorem~\ref{theorem_skip}.
Let us consider a situation where the data fidelity term $\Phi=\|\bm{A}\cdot-\bm{y}\|^2$ wants to be minimized under the prior induced by the denoiser, where $\bm{A}\in\mathbb{R}^{d\times d_0}$ and $\bm{y}\in\mathbb{R}^{d}$ are a given matrix and a given vector, respectively.
The well-known soft thresholding operator can be used as a denoiser \cite{donoho1995_softshrink}, since it corresponds to the proximity operator of the convex regularizer $\|\cdot\|_1$.
However, it tends to cause underestimation, and therefore many nonconvex regularizers have been proposed \cite{Fan2001_SCAD_shrink, chartrand2007_lp, candes2007_reweighted_l1, Zhang2010_MCP_shrink, yao2017_family_nonconvex, yukawa2023_limes}.
Inspired by those previous works, the MoL-Grad denoiser is constructed as the proximity operator (in the sense of Definition~\ref{def:sprox}) of a weakly convex regularizer, which enhances the sparsity-promoting property.
This means that plugging MoL-Grad into a splitting algorithm leads to minimization of a cost function involving its associated weakly convex regularizer under appropriate conditions.
Nevertheless, to preserve the overall convexity of the cost function, the data fidelity term $\Phi$ must be strongly convex.
In a typical image restoration problem, however, the matrix $\bm{A}$ is often underdetermined (i.e., $d < d_0$) or ill-posed (i.e., $\bm{A}^\mathsf{T}\bm{A}$ has negligible eigenvalues).
In such a situation, the MoL-Grad denoiser cannot be directly applied.
To address this issue, we propose to restrict the effect of the enhancement to a subspace $\mathcal{M}$ on which convexity of the whole cost can be guaranteed.
This restriction strategy is embodied by adding a squared norm penalty $\|P_{\mathcal{M}^\perp}\cdot\|^2$ on the orthogonal complement $\mathcal{M}^\perp$.
The algorithm is summarized in Algorithm~\ref{alg:primal_dual}, where $D$ is the MoL-Grad denoiser.

\begin{algorithm}[t]
\caption{Primal-Dual Algorithm for Image Restoration}
\label{alg:primal_dual}
\begin{algorithmic}[1]
\State \textbf{Set}: $\bm{v}_0, \bm{u}_0\in\mathbb{R}^{d_0}, \sigma, \tau, \mu >0$.
\For{$k = 0, 1, 2, 3, \cdots $}
        \State $\tilde{\bm{u}}_{k+1}:= \bm{u}_k + \sigma \bm{v}_k$
        \State $\bm{u}_{k+1}:= \tilde{\bm{u}}_{k+1} - \sigma D (\tilde{\bm{u}}_{k+1}/(\sigma+1))$
        \State $\tilde{\bm{v}}_{k+1}:=  \tau \nabla (\mu \Phi+(1/2)\|P_{\mathcal{M}^\perp}\cdot\|^2)(\bm{v}_k)$
        \State $\bm{v}_{k+1}:=(1+\tau) \bm{v}_k -\tilde{\bm{v}}_k-\tau (2\bm{u}_{k+1}-\bm{u}_k)$
\EndFor

\State \textbf{return} $\bm{v}_{k+1}$  
\end{algorithmic}
\end{algorithm}
\begin{prop}{}
\label{prop:conv_ana}
Let $D:\mathbb{R}^{d_0}\to\mathbb{R}^{d_0}:\bm{v}\mapsto\nabla\psi(\bm{v})$ be a MoL-Grad denoiser for some $\psi\in\Gamma_0(\mathbb{R}^{d_0})$, i.e., $D$ satisfies condition $\clubsuit$ for some $\beta\in(0,1)$.
In this case, $D$ is the proximity operator of the ``implicit'' regularizer $\varphi:=\psi^*-(1/2)\|\cdot\|^2$ (see Fact~\ref{fact:def_molgrad}).
Let $\Phi\in\Gamma_0(\mathbb{R}^{d_0})$ be a smooth function, which is $(1/\mu)$-strongly convex for some $\mu\in\mathbb{R}_{++}$ over a subspace $\mathcal{M}$, so that the function
\begin{align}
    \Gamma_0(\mathbb{R}^{d_0})\ni h: \bm{v}\mapsto \mu\Phi(\bm{v})-(1/2)\|P_{\mathcal{M}}\bm{v}\|^2 \label{eq:def_h}
\end{align}
is $\kappa$-smooth for some $\kappa \in \mathbb{R}_{++}$.
Here, $P_\mathcal{M}$ and $P_{\mathcal{M}^\perp}$ denote the orthogonal projections onto $\mathcal{M}$ and its orthogonal complement $\mathcal{M}^\perp$, respectively.
Set the step sizes $\sigma, \tau>0$ such that
\vspace{5pt}\\
$
\text{\normalfont (i)}  \hspace{10pt} \displaystyle \sigma \leq \frac{ \beta}{1-\beta}
$
\vspace{5pt}\\
$
\text{\normalfont(ii)} \hspace{10pt} \displaystyle \tau\left( \sigma + \frac{\kappa}{2} \right) < 1.
$
\vspace{5pt}\\
Then, for arbitrary $(\bm{v}_0, \bm{u}_0) \in \mathbb{R}^{d_0}\times\mathbb{R}^{d_0}$, the sequences $(\bm{v}_k)_{k\in\mathbb{N}}\subset\mathbb{R}^{d_0}$ and $(\bm{u}_k)_{k\in\mathbb{N}}\subset\mathbb{R}^{d_0}$ generated by Algorithm \ref{alg:primal_dual} converge to solutions $\check{\bm{v}}\in\mathbb{R}^{d_0}$ and $\check{\bm{u}}\in\mathbb{R}^{d_0}$ of the following primal and dual problems, respectively:
\begin{align}
    &\min_{\bm{v}\in\mathbb{R}^{d_0}} \mu\Phi(\bm{v})+(\sigma+1)\varphi(\bm{v})+\frac{1}{2}\|P_{\mathcal{M}^\perp}(\bm{v})\|^2,\label{eq:primal problem}\\
    &\min_{\bm{u}\in\mathbb{R}^{d_0}} h^*(\bm{u})+\Big[(\sigma+1)\varphi+\frac{1}{2}\|\cdot\|^2\Big]^*(\bm{u}),
\end{align}
provided that such solutions exist.
\vspace{0.5em}\par\noindent {\itshape Proof.}\hspace{7pt}\normalfont
The assesrtion can be verified straightforwardly
by setting $L := \bm{I}$, $T:=D$, $\rho:=1$, $\hat{f}:=h$, and $f:=h+(1/2)\|\cdot\|^2=\mu\Phi+(1/2)\|P_{\mathcal{M}^\perp}\cdot\|^2$ in \oldcite[Algorithm~1]{yukawa2025_molgrad} and by applying \oldcite[Theorem~3]{yukawa2025_molgrad}.
\qed
\end{prop}
In \eqref{eq:primal problem}, the addition of the extra penalty $(1/2)\|P_{\mathcal{M}^\perp}\cdot\|^2$ allows to enhance the strength of the implicit regularizer $\varphi$ induced by the denoiser $D$.
Specifically, this modification relaxes the upper bound of $\sigma$, so that overall convexity of the cost function in \eqref{eq:primal problem} can be guaranteed for larger values of $\sigma$.
This technique is inspired by the approach discussed in \cite{suzuki_external}.\footnote{
We consider a slightly more general case
in which $\bm{A}^\mathsf{T}\bm{A}$ is full-rank but involves
small eigenvalues;
\cite{suzuki_external} focuses on the underdetermined case.}
\begin{example}{Selection of $\mathcal{M}$ for the quadratic data fidelity setting}
Let us consider a case of
\begin{equation}
    \Phi(\bm{v})=(1/2)\|\bm{A}\bm{v}-\bm{y}\|^2, \label{eq:def_f_as_normsq}
\end{equation}
where $\bm{A}\in\mathbb{R}^{d\times d_0}$ and $\bm{y}\in\mathbb{R}^{d_0}$ with $\lambda_1(\bm{A}^\mathsf{T}\bm{A})\geq \lambda_2(\bm{A}^\mathsf{T}\bm{A}) \geq \cdots \geq \lambda_i(\bm{A}^\mathsf{T}\bm{A})\geq 1/\mu \geq \lambda_{i+1}(\bm{A}^\mathsf{T}\bm{A}) \geq \cdots \geq \lambda_{d_0}(\bm{A}^\mathsf{T}\bm{A}) (\geq 0)$.
In this case, $\mathcal{M}$ is the dominant eigenspace corresponding to the largest $i$ eigenvalues under which the inclusion condition in \eqref{eq:def_h} is satisfied.
\end{example}
%
%

%

\begin{remark}{Guideline of parameter setting}
\label{remark:LD}
Algorithm~\ref{alg:primal_dual} has three parameters $\sigma$, $\mu$, and $\tau$, to be designed.
We choose $\tau=\gamma (\sigma+(\kappa/2))^{-1}$ where $\gamma\in(0, 1)$ so that condition (ii) of Proposition \ref{prop:conv_ana} is satisfied.
In practice, one can fix the parameter $\gamma$ to some constant, say $\gamma := 0.8$, because $\tau$ only affects the convergence speed.
%
Viewing \eqref{eq:primal problem},
on the other hand, $\sigma$ and $\mu$
determine the weights of the implicit regularizer $\varphi$ and the
data fidelity term $\Phi$, respectively.
As a rule of thumb, better performance tends to be obtained by assigning
a larger weight to $\varphi$ compared to the other terms.
If the Lipschitz constant $L_D$ is available, a recommended strategy
would thus be to set $\sigma$ to the largest possible value $\sigma = \beta /
(1 - \beta)$ according to condition (i) of Proposition \ref{prop:conv_ana},
where $\beta=L_D^{-1}$.
The remaining parameter $\mu$ needs to be tuned manually.

In practice, the Lipschitz constant can be estimated by
$(L_D \geq )$~$\hat{L}_D:=\max_{\bm{x}
\in X} \{\| \mathrm{J}_D(\bm{x}) \|\}$,
where $X$ is a set of noisy input images
(see \cite{pesquet2021_MMO},\cite{hurault2022proximal}).\footnote{
In Theorem~\ref{theorem_mol_NN}(b), the Lipschitz constant of $D = R_1$
 can be estimated using Lemma~\ref{memo_lip2}, based on
 $L_{\widetilde{T}_{n}}$, $L_{R_{n}}$ and $\|\bm{W}_n^\mathsf{T}\|$.
See Appendix \ref{app:proof_of_theorem1} for the definitions
of the operators $\widetilde{T}_{n}$ and $R_{n}$.
As can be seen from the discussions of the appendix
as well as Lemma \ref{memo_lip2},
$L_D$ can be estimated from $\widetilde{T}_n$,
since $R_{n+1}$ is defined recursively with
$R_N=\bm{W}_N^\mathsf{T}\widetilde{T}_N$.
Given that $\widetilde{T}_n$ has a general neural network structure,
the method proposed in \cite{combettes_lip} can be applied to
estimation of tight Lipschitz constants.
However, even if a tight Lipschitz constant of $\widetilde{T}_n$ is
provided, the resulting $L_D$ computed via Lemma~\ref{memo_lip2} tends
to be large, as it is obtained via a product of multiple constants.
}
As this practical estimate gives a lower bound of the true $L_D$,
a possible way of setting $\sigma$ is to use $\hat{L}_D+\epsilon$
for some small $\epsilon>0$.
This conservative strategy however yields small $\sigma$,
resulting in suboptimal performance.
We therefore use $\hat{L}_D - \tilde{\epsilon}$
for some $\tilde{\epsilon}>0$ which makes the corresponding $\sigma$
larger in our simulations, although the convergence conditions are
violated in this case.

\end{remark}

The forward-backward splitting type algorithm employing the MoL-Grad denoiser has the range of the regularization parameter strictly bounded from below
(as well as from above) by a constant depending on the Lipschitz constant
$L_D$ \oldcite[Remark 3]{yukawa2025_molgrad}.
This motivates us to adopt the primal-dual type algorithm that is free from such a restriction.
\section{Simulation Results}
\label{sec:simulation}
We present a simulation studies of the proposed denoiser from three perspectives.
We first analyze the gain coming from the Lipschitz-constraint-free
nature of the proposed denoiser.
We then show that the proposed denoiser enjoys significantly shorter training time.
We finally show that it gives the state-of-the-art performance
in image deblurring among the existing PnP methods that have convergence guarantees.
%
%

\subsection{Denoiser Analysis}
\label{subsec:result_denoise}
\begin{figure}[t]
    \centering
    \includegraphics[width=0.85\linewidth]{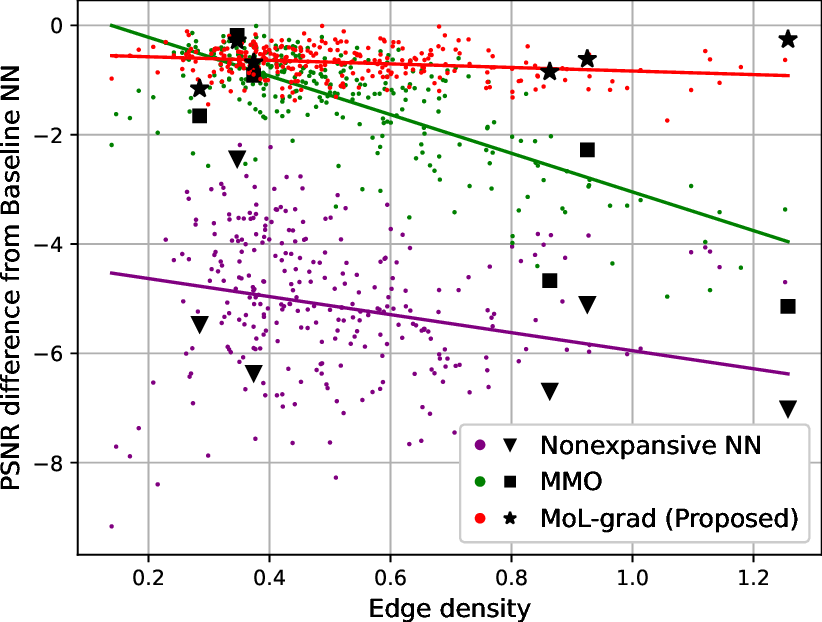}
    \vspace{5pt}
    {\captionsetup[subfloat]{labelformat=empty}
    \subfloat[(i) bag]{\includegraphics[width=0.26\linewidth]{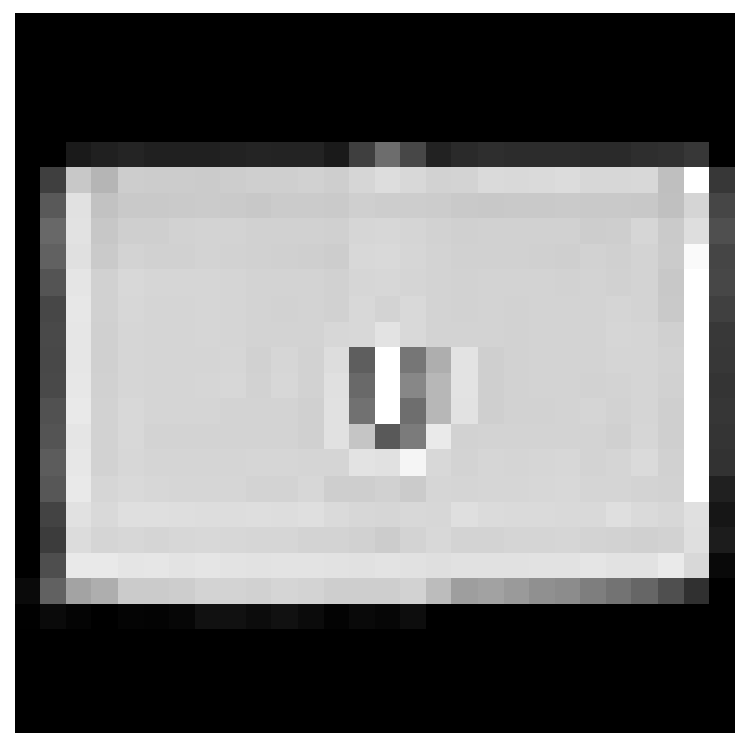}}
    \hfill
    \subfloat[(ii) coat]{\includegraphics[width=0.26\linewidth]{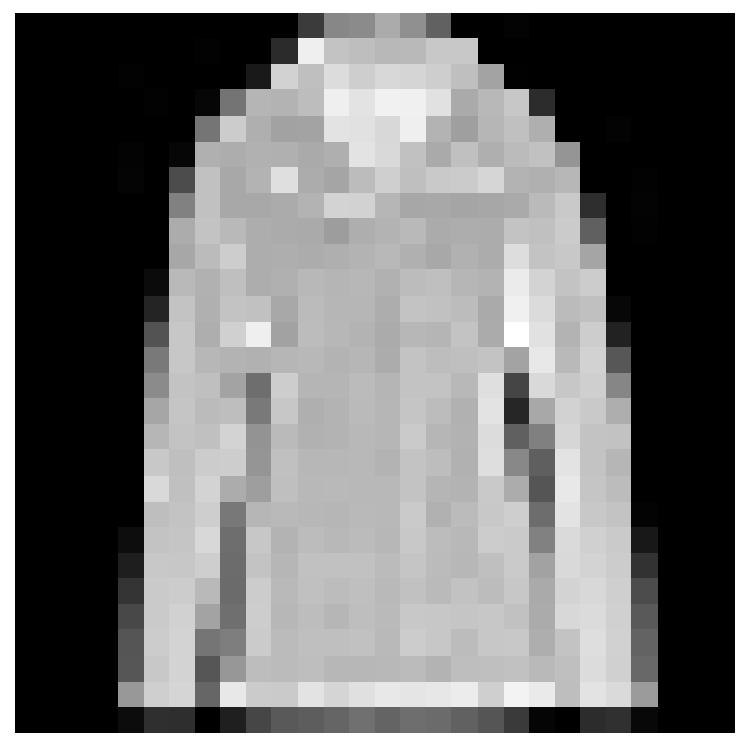}}
    \hfill
    \subfloat[(iii) boot]{\includegraphics[width=0.26\linewidth]{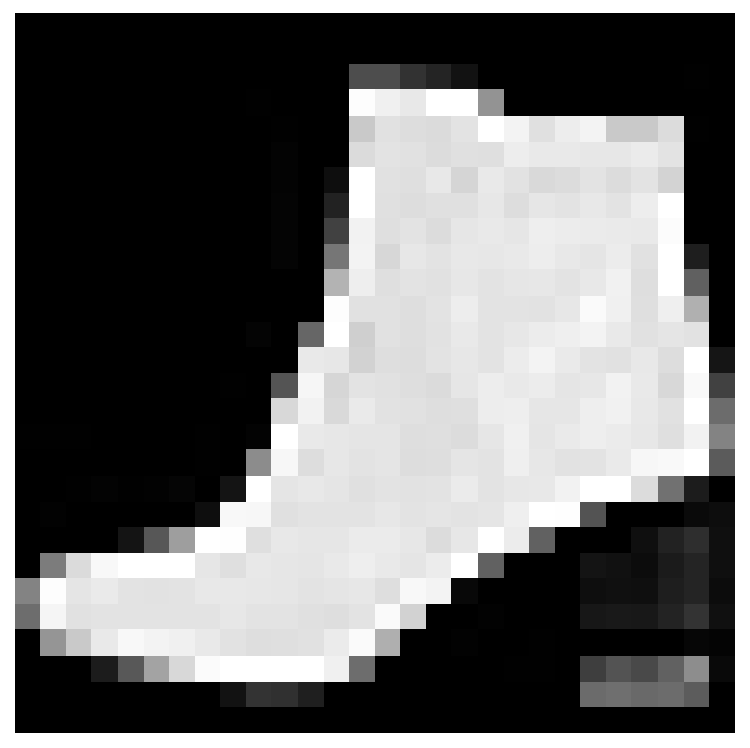}}

    \vspace{2pt}

    \subfloat[(iv) shoe]{\includegraphics[width=0.26\linewidth]{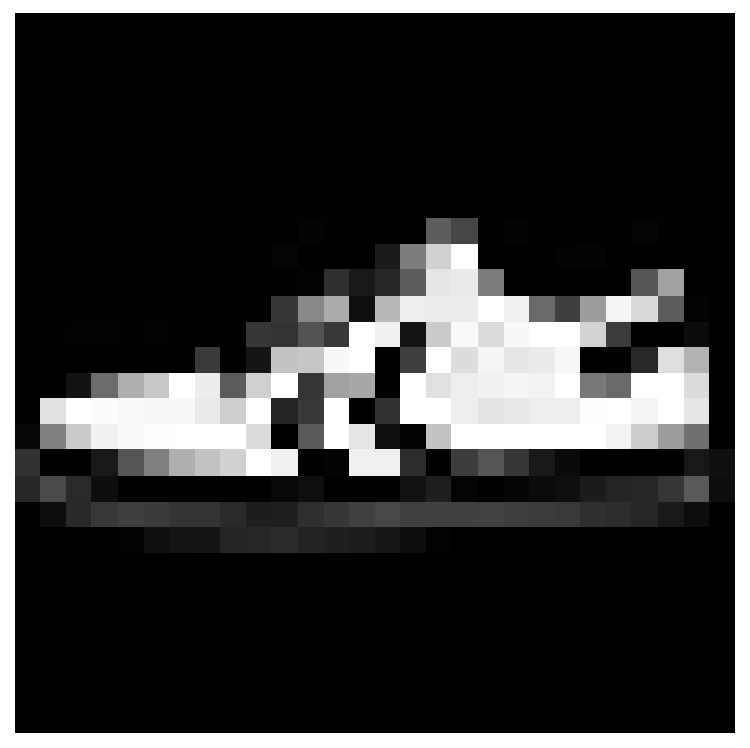}}
    \hfill
    \subfloat[(v) top]{\includegraphics[width=0.26\linewidth]{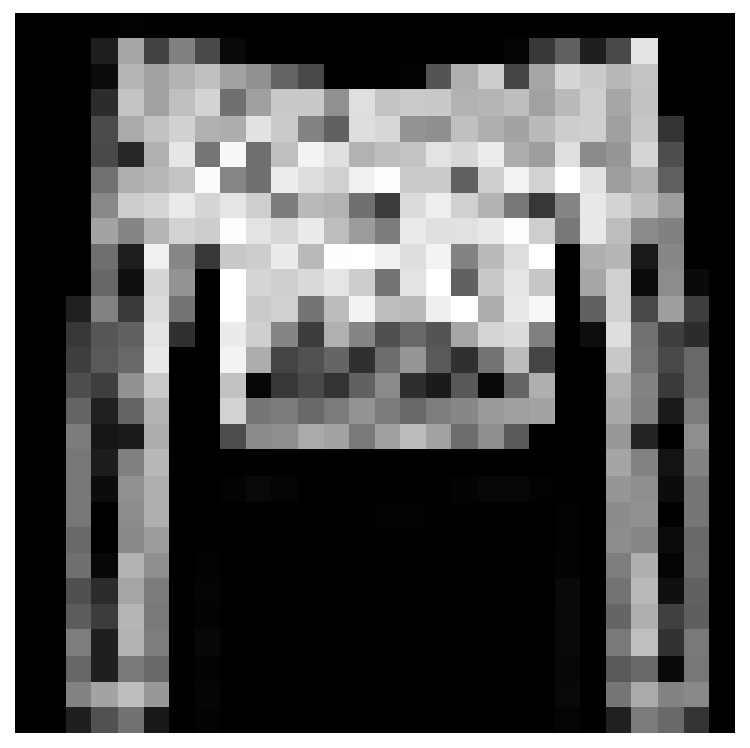}}
    \hfill
    \subfloat[(vi) sandal]{\includegraphics[width=0.26\linewidth]{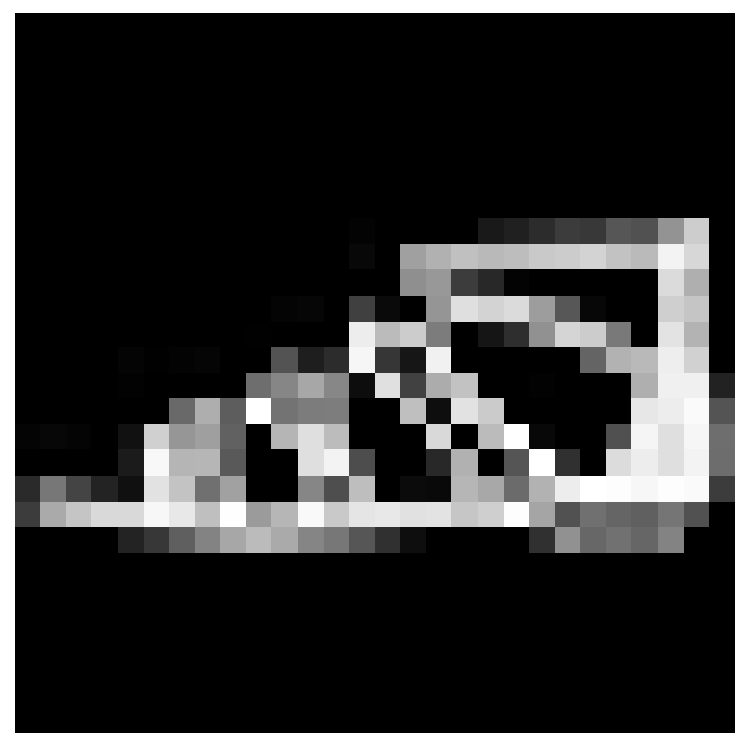}}
    }

    \caption{PSNR Difference from baseline NN vs. Edge density. The line represents the least-squares linear fit. The black markers ({\large $\blacktriangledown$}, \hspace{2pt}\rule{1.5ex}{1.5ex}\hspace{2pt}, $\bigstar$) represent, from left to right, the images shown below in (i) to (vi), respectively.}
    \label{fig:edge_psnr}
\end{figure}
\begin{table}[tpb]
\centering
\caption{Denoising performance for the images presented shown in Fig.~\ref{fig:edge_psnr}}
{\setlength{\tabcolsep}{3pt}

\subfloat[Images with low edge density (PSNR in dB)]{
  {\fontsize{10pt}{12pt}\selectfont
  \begin{tabular}{c|ccc}
    & (i) & (ii) & (iii) \rule{0pt}{2ex}\\
    \hline
    \rule{0pt}{2ex}
    Edge density & 0.29 & 0.35 & 0.37\\
    \hline\hline
    \rule{0pt}{2ex}
    Nonexpansive NN            & 27.61 & 28.40  & 25.88\\
    MMO                  & 31.43 & \textbf{30.66} & 31.33\\
    MoL-Grad (Proposed)  & \textbf{31.92} & 30.57 & \textbf{31.57}\\
  \end{tabular}
  }
  \label{tab:psnr_simple}
}
\hfill
\subfloat[images with high edge density (PSNR in dB)]{
  {\fontsize{10pt}{12pt}\selectfont
  \begin{tabular}{c|ccc}
    & (iv) & (v) & (vi) \rule{0pt}{2ex}\\
    \hline
    \rule{0pt}{2ex}
    Edge density & 0.86 & 0.93 & 1.26\\
    \hline\hline
    \rule{0pt}{2ex}
    Nonexpansive NN            & 25.17 & 23.70 & 24.05 \\
    MMO                  & 27.20 & 26.53 & 25.87 \\
    MoL-Grad (Proposed)  & \textbf{31.02} & \textbf{28.19} & \textbf{30.78} \\
  \end{tabular}
  }
  \label{tab:psnr_complex}
}
}
\label{table:psnr_denoise}
\end{table}
\begin{figure}[tb]
  \centering

  \subfloat[True]{%
    \begin{minipage}{0.48\linewidth}
      \centering
      \includegraphics[width=\linewidth]{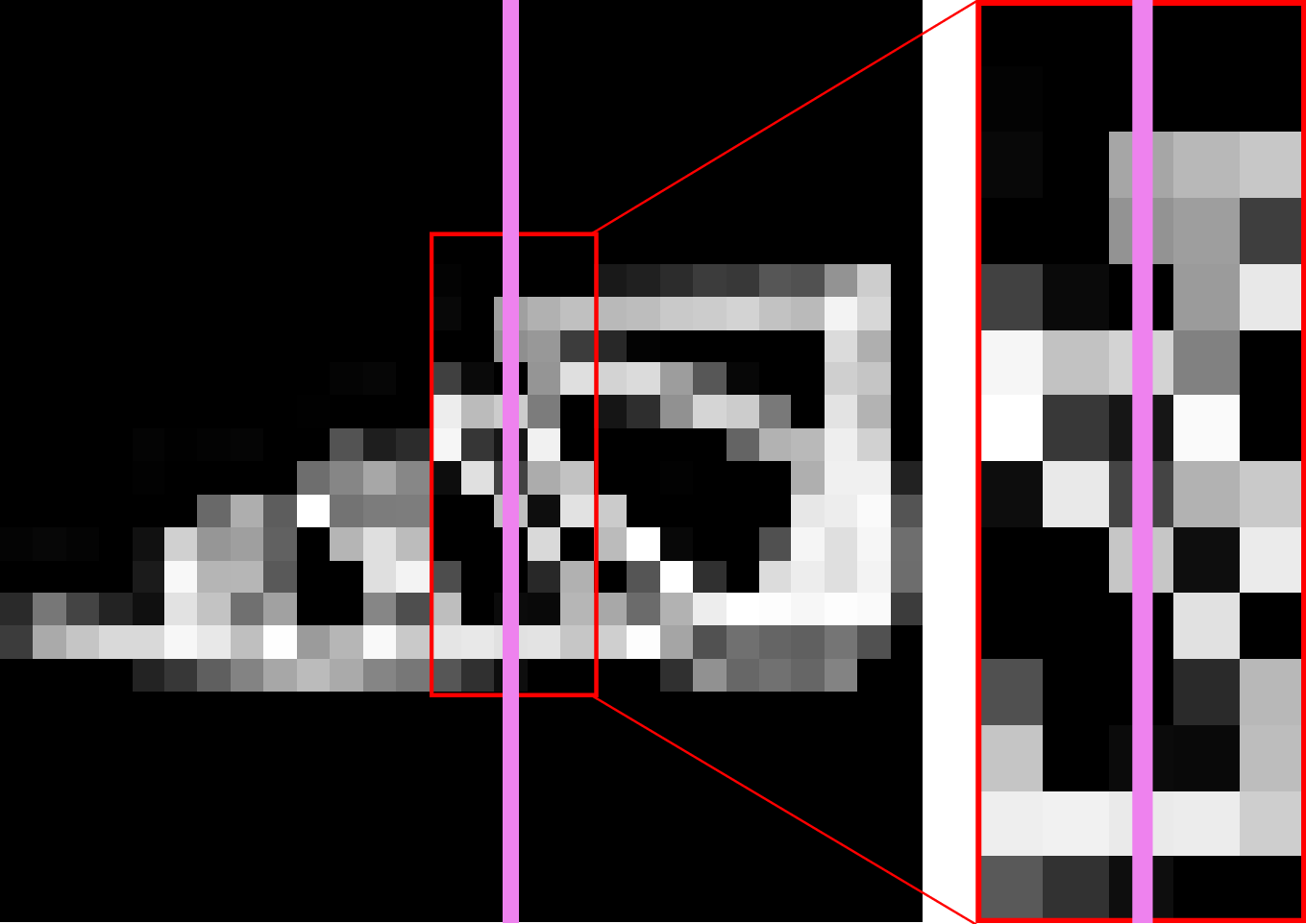}\\
      \hspace{3pt}
      \includegraphics[width=\linewidth]{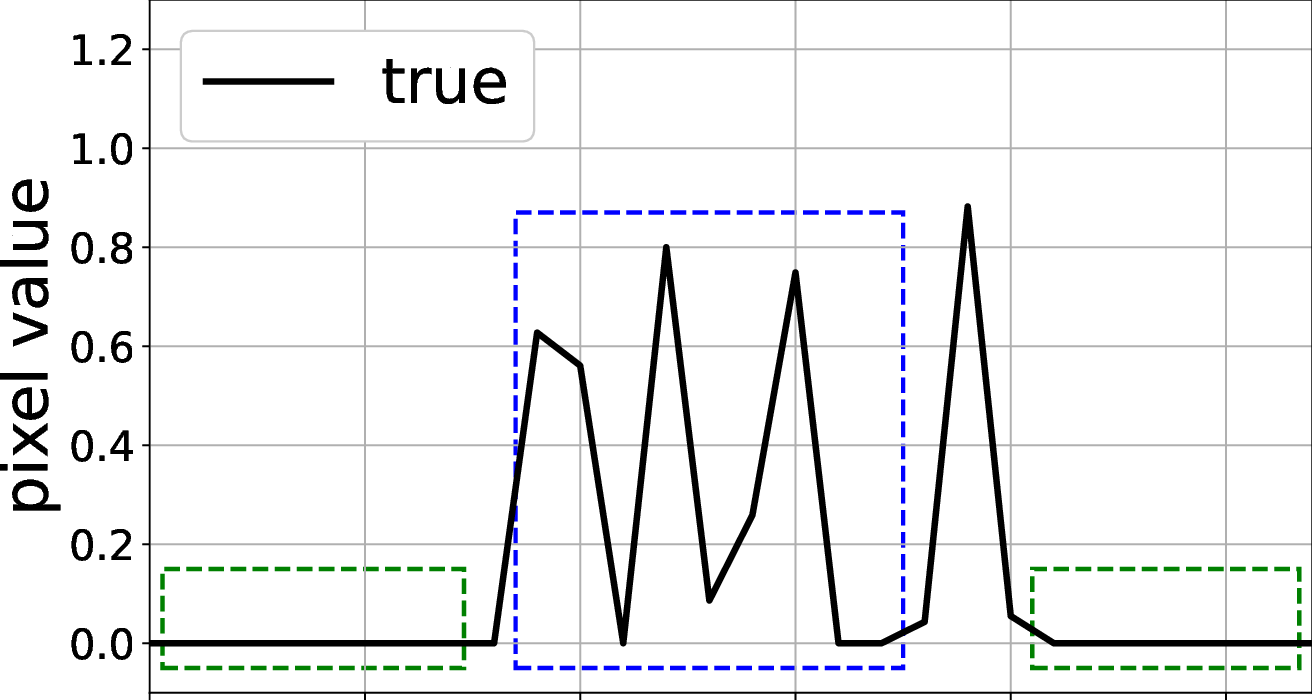}
    \end{minipage}
  }\hfill
  \subfloat[Nonexpansive NN\\ \hspace{10pt}(PSNR $=24.05$)]{%
    \begin{minipage}{0.48\linewidth}
      \centering
      \includegraphics[width=\linewidth]{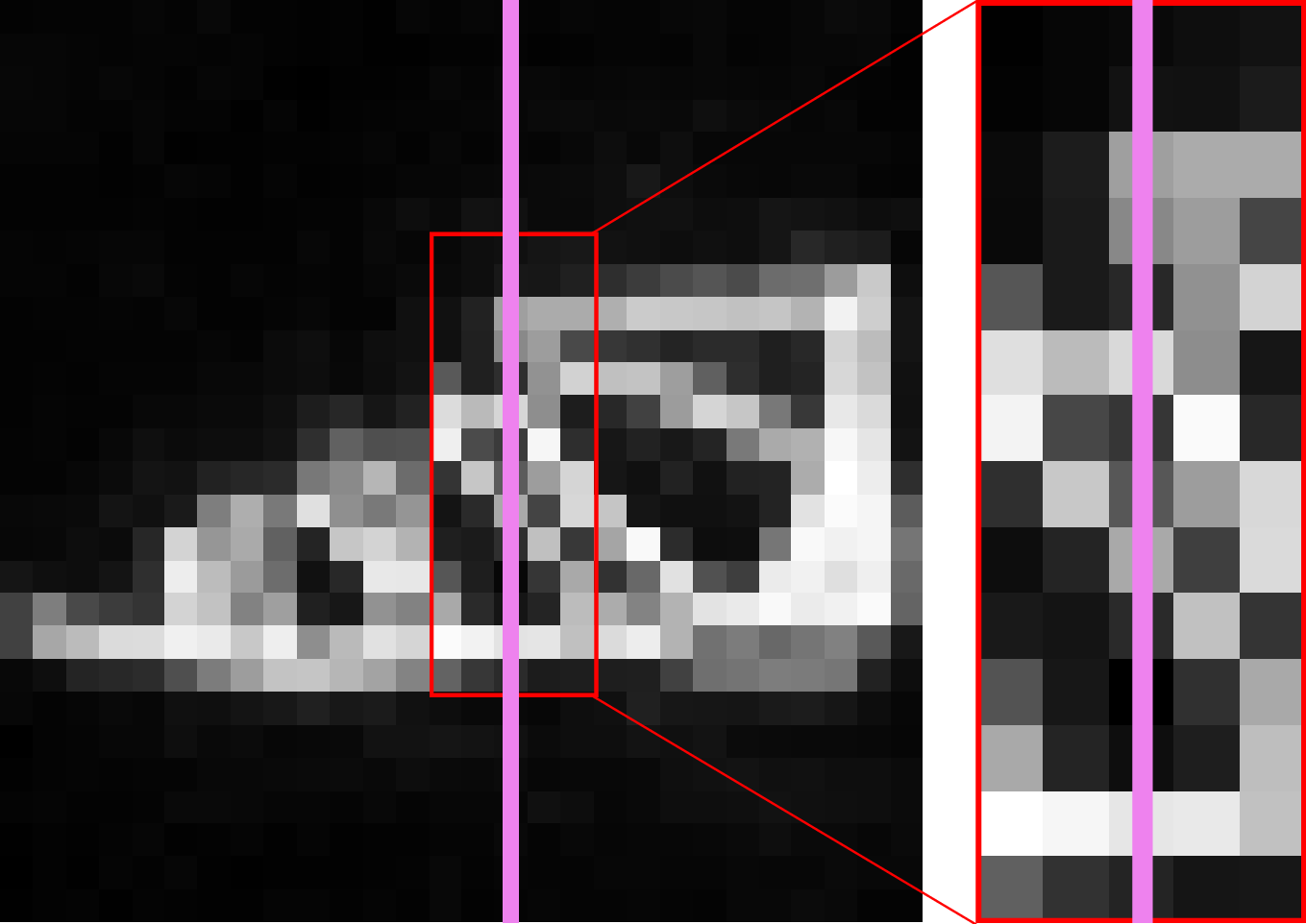}\\
      \hspace{3pt}
      \includegraphics[width=\linewidth]{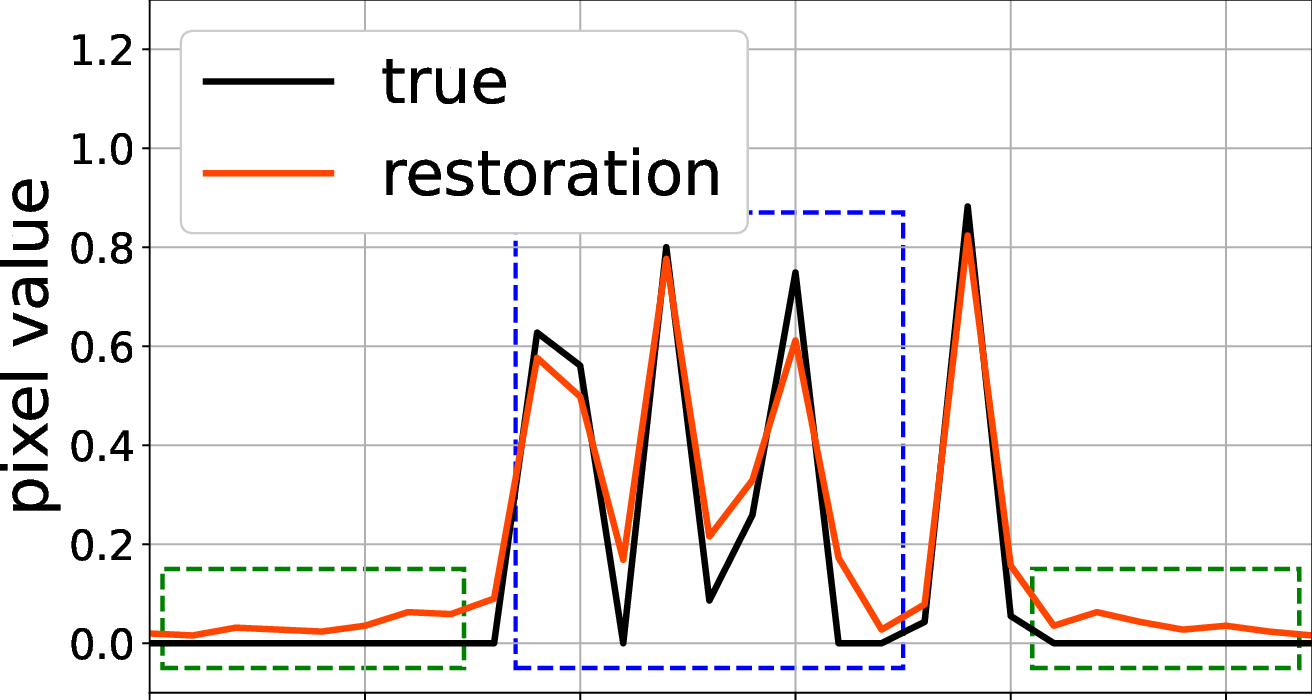}
    \end{minipage}
  }


  \subfloat[MMO (PSNR $=25.87$)]{%
    \begin{minipage}{0.48\linewidth}
      \centering
      \includegraphics[width=\linewidth]{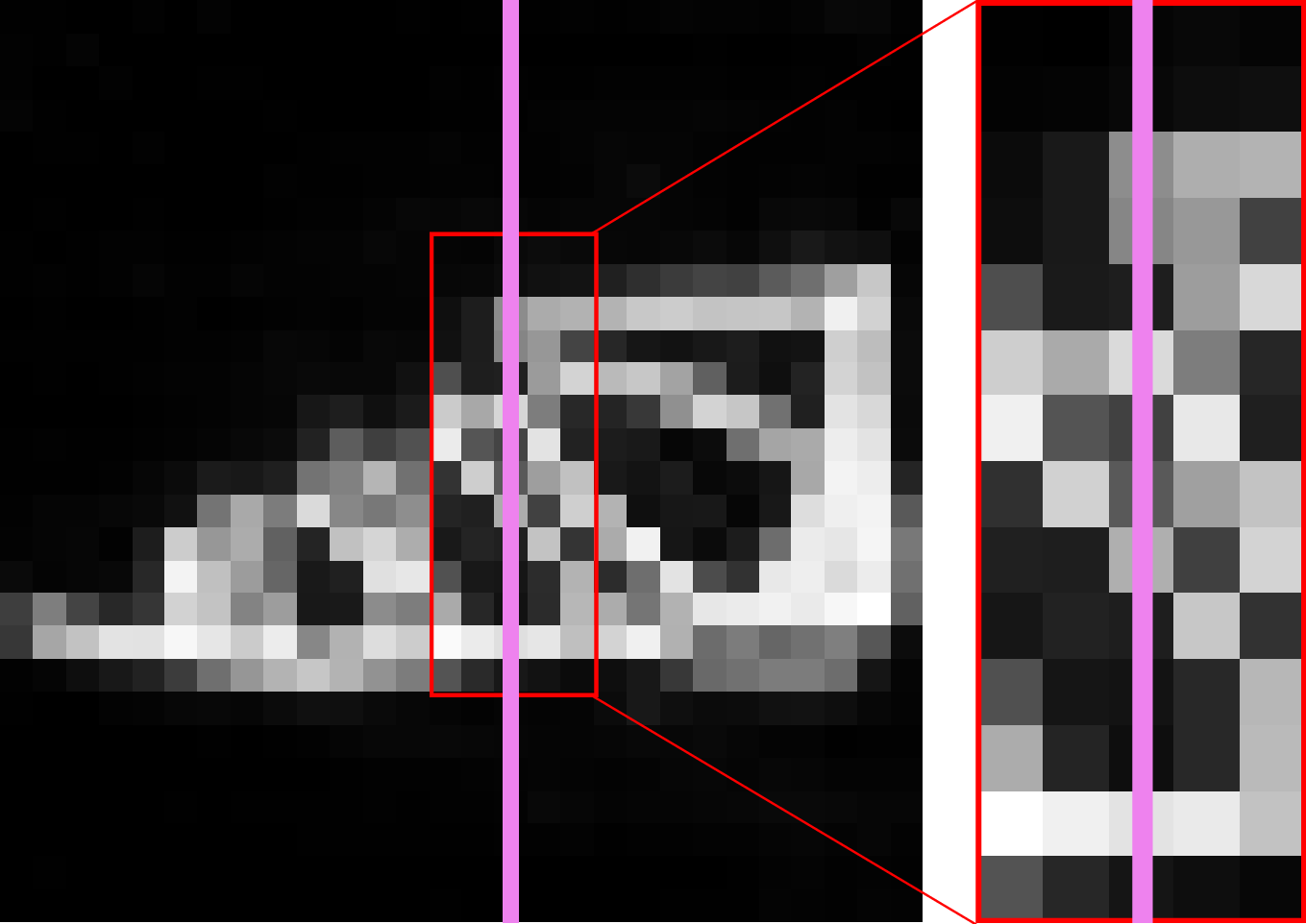}\\
      \hspace{3pt}
      \includegraphics[width=\linewidth]{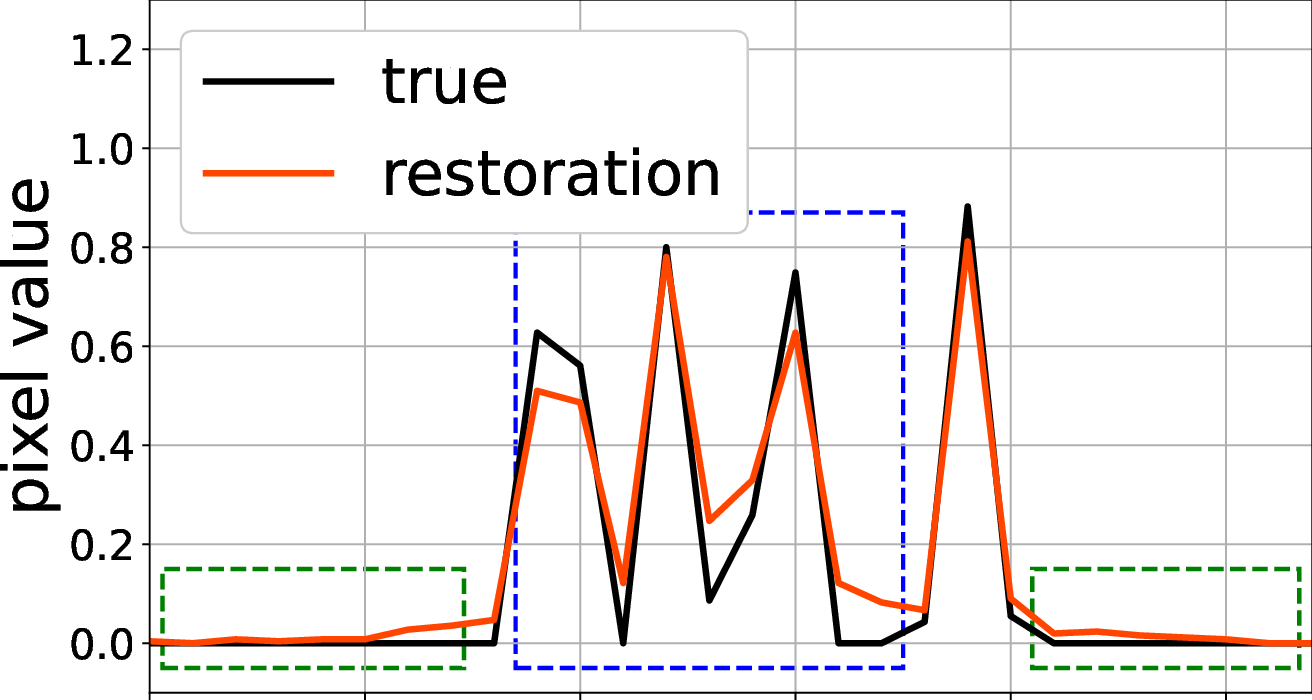}
    \end{minipage}
  }\hfill
  \subfloat[MoL\!-Grad (Proposed)\\ \hspace{10pt}(PSNR $=30.77$)]{%
    \begin{minipage}{0.48\linewidth}
      \centering
      \includegraphics[width=\linewidth]{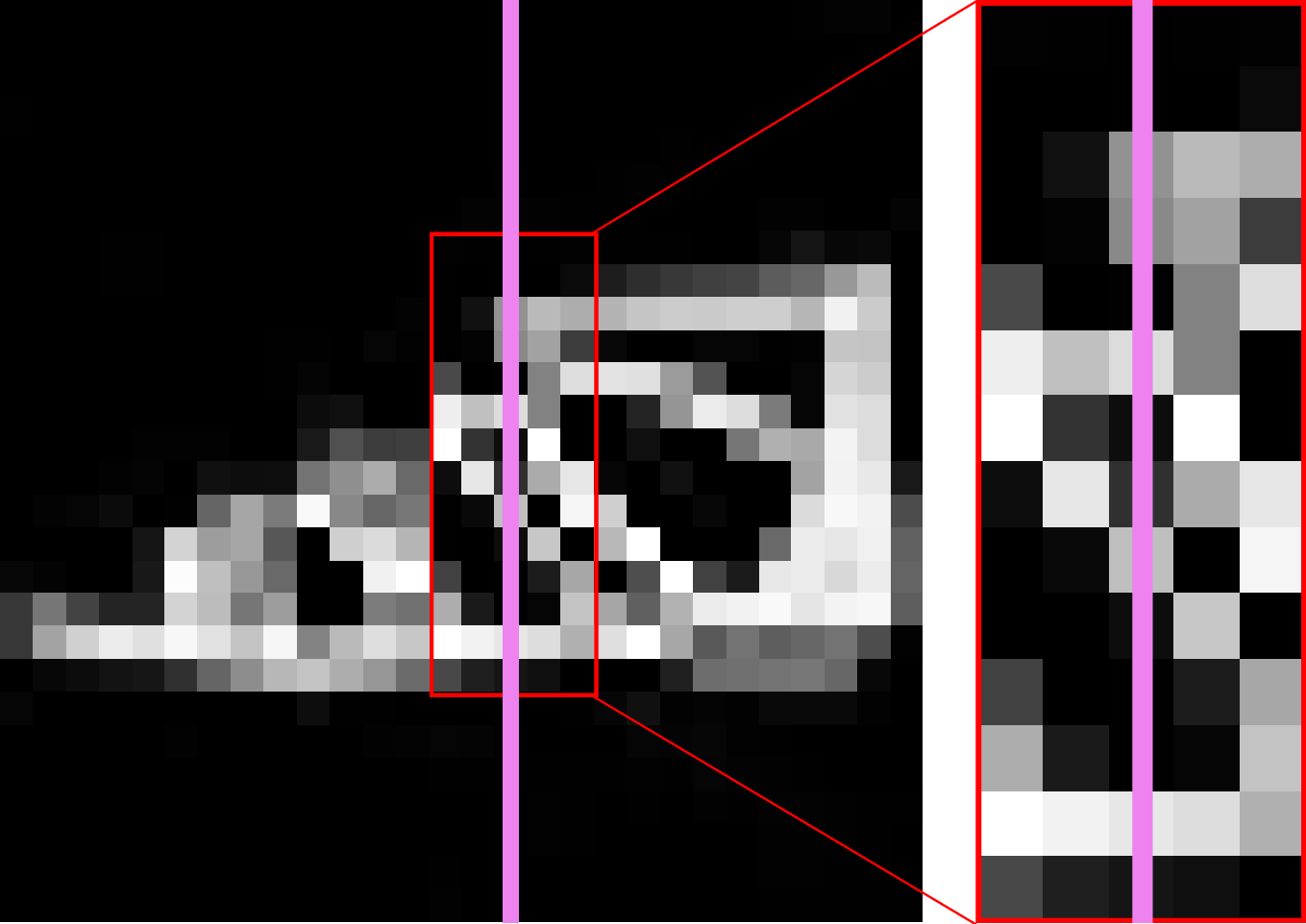}\\
      \hspace{3pt}
      \includegraphics[width=\linewidth]{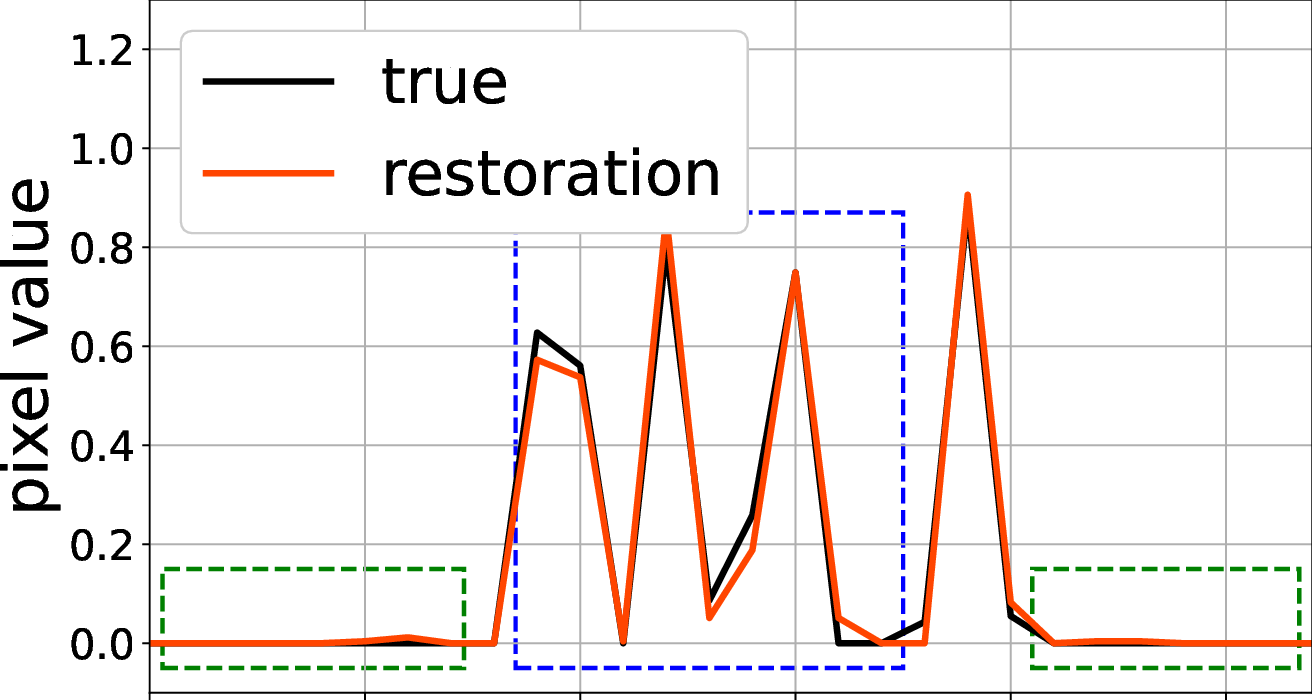}
    \end{minipage}
  }

  \caption{Denoising results of ``sandal''. The graphs show
 the pixel values along the pink lines (PSNR in dB).}
  \label{fig:denoise_sandal}
\end{figure}
{\it Can the proposed MoL-Grad denoiser (which has no Lipschitz constant
bound) outperform the nonexpansive ones
 (see Remark \ref{remark:nonexpansive_limitation})?}
To address this question, simulations are set up as below.
For training, $6000$ images are used from the Fashion-MNIST dataset
\cite{Xiao2017_fMNIST} with the noise level $\sigma_\varepsilon = 0.05$.
The learning rate is set to $1.0 \times 10^{-4}$ for the first 800 training epochs and to $2.5\times10^{-6}$ for the remaining 200 epochs.
The proposed denoiser is compared to two nonexpansive denoisers:
nonexpansive neural network (Nonexpansive NN) and
maximally monotone operator (MMO)
\cite{pesquet2021_MMO}.\footnote{Nonexpansive NN is trained in a way that  $\hat{L}_D \leq 1$ is satisfied.}
%
%
For the proposed method,
the network is trained in a way the assumptions of Theorem
\ref{theorem_mol_NN} are strictly satisfied,
whereas those for Nonexpansive NN and MMO are trained
subject to the relaxed Lipschitz-constant constraints
(cd. Remark~\ref{remark:LD}).
We define ``edge density'' $\Omega(\bm{x})$ to quantify the ``structural complexity'' of an image as follows:
\begin{align}
    \Omega(\bm{x}) &:= {\rm supp}(\bm{x})^{-1}
\sum_{\delta\in \{{\rm h,v}\}}
\sum_{i} H(|(\bm{D}_\delta\bm{x})_i|-t_{\mathrm{th}}),
\end{align}
where $H$ denotes the Heaviside step function,
$\bm{D}_{\rm h}$ and $ \bm{D}_{\rm v}$ are the first-order difference
operators in the horizontal and vertical directions, respectively, and
$t_{\mathrm{th}} = 0.05$ is the threshold for edge magnitude.
Fig.~\ref{fig:edge_psnr} showcases the relationship between edge density
and the denoising performance for 300 images (different from
those used for training) in the Fashion-MNIST dataset,
where the vertical axis shows the difference from the unconstrained neural network (baseline NN).
To facilitate the analysis, Table~\ref{table:psnr_denoise}
highlights the denoising performance for selected images
from two categories:
images of high edge density (complex images)
and images of low edge density (simple images).
It can be seen that the proposed denoiser yields significant gains
for the complex images.
Fig.~\ref{fig:denoise_sandal} shows the denoised images of
``sandal'' (a specific example of the complex images in the dataset).
To enhance the resolution of the analysis,
let us focus on the specific part of the image along the line in pink
color.
The graphs below the images plot the values of the pixels on the line.
Inspecting the green boxes in the graphs,
the proposed denoiser suppresses the noise appearing
in the black part of the image successfully,
while Nonexpansive NN completely fails and
MMO still has some noise remaining in the restored image.
Turning our attention to the blue boxes in the graphs,
the proposed denoiser restores the fine structures of the sandal image
much better than the other denoisers.
To sum up, thanks to the expansiveness coming from the
Lipschitz-constraint-free nature,
the proposed denoiser effectively preserves the edges
while achieving superior denoising performance simultaneously.

%
\subsection{Training Speed}
\label{subsec:training_speed}
%
%
%
%
\begin{algorithm}[t]
\caption{Training algorithm used in Section~\ref{subsec:training_speed}}
\label{alg:training}
\begin{algorithmic}[1]
\State Let $K\in\mathbb{N}$ be the number of epoch size.
\For{$k = 1, 2, 3, \cdots , K$}
    \State \hspace{-1pt}{\bf Forward step:} $\bm{x}^*=D_{\bm{\theta}_k} (\bm{x}+\bm{\varepsilon})$
    \State \hspace{-1pt}{\bf Loss step:} Calculate $L(\bm{x}, \bm{x}^*, \bm{\theta}_k)$
    \State \hspace{-1pt}{\bf Backward Step:} $\bm{\theta}_{k+1} = \text{Adam}(\nabla_{\theta_k} L(\bm{x}, \bm{x}^*, \bm{\theta}_k), \bm{\theta}_k)$
\EndFor
\State \textbf{return} $\bm{\theta}_{k+1}$
\end{algorithmic}
\end{algorithm}
{\setlength{\tabcolsep}{1pt}
\begin{table}[t]
  \caption{Computation time for each training step of one epoch (millisecond).}
  \centering
  {\fontsize{9pt}{11pt}\selectfont
  \begin{tabular}{c|ccc|cc}
    \hline
    \rule{0pt}{2.5ex}
    & \hspace{1pt}Forward & Loss & Backward & {\bf Total}\\
    \hline\hline
    \rule{0pt}{2.5ex}
    Baseline NN \hspace{1pt} & $6.16$ & $3.00\!\times\!10^{-2}$ & $1.00\!\times\! 10^1$ & $1.62\!\times\!10^1$ \\
    \hline
    \rule{0pt}{2.5ex}
    MMO & $6.39$ & $5.13\!\times\! 10^{1}$ & $3.35\!\times\!10^{2}$ & $3.93\!\times\!10^{2}$\\
    Prox-GS & $10.73\hspace{3pt}$ & $3.62\!\times\!10^{3}$ & $6.95\!\times\!10^{2}$ & $4.33\!\times\!10^{3}$ \\
    MoL-Grad & $\bm{5.46}$ & $\bm{8.50\!\times\! 10^{-1}}$ & $\bm{1.60\!\times\!10^{2}}$ & $\bm{1.66\!\times\!10^{2}}$\\
    \hline
  \end{tabular}
  }
  \label{tab:training_speed}
\end{table}
}
\begin{figure}[tbp]
    \centering
    \subfloat[Case of $128\times128$ pixels]{
        \includegraphics[width=0.8\linewidth]{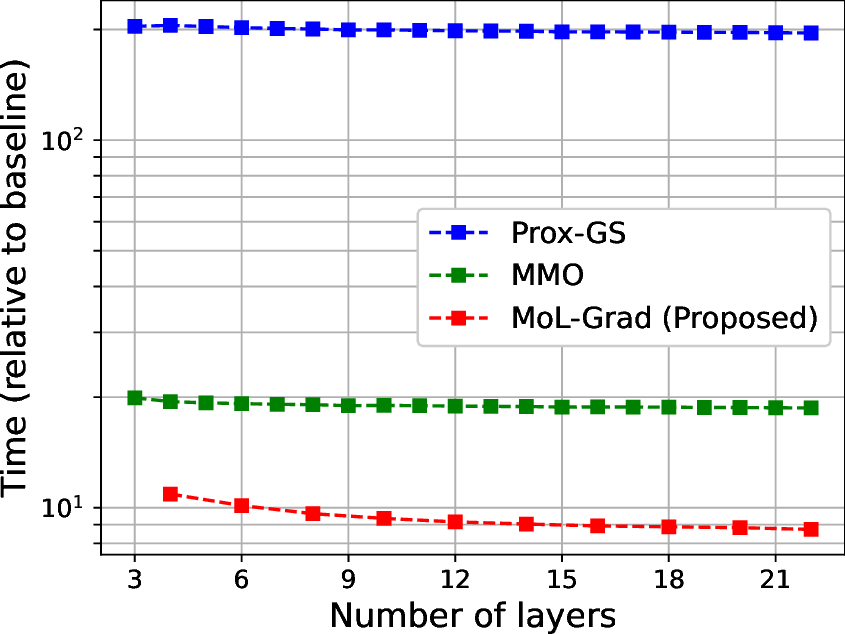}
        \label{fig:time_vs_layers}
    }

    \vspace{6pt}
    \subfloat[Case of $4$ layer neural network]{
        \includegraphics[width=0.8\linewidth]{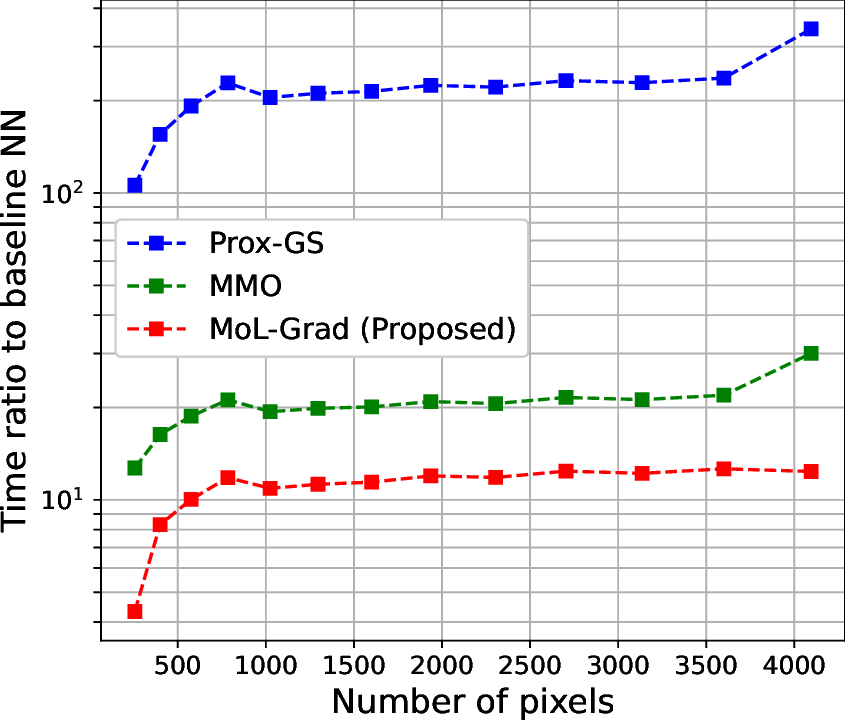}
        \label{fig:time_vs_input}
    }

    \caption{Comparisons in training time relative to the baseline NN under various conditions.
The number of layers for the proposed method is an even number because
 of the weight-tying structure.}
    \label{fig:training_time}
\end{figure}

The Lipschitz-constraint-free nature of the proposed denoiser
permit its training without time-consuming power iteration
which is required for training those denoisers which has convergence
guarantee.
We show that this yields significant reduction of the training time
compared to MMO \cite{pesquet2021_MMO} and Prox-GS \cite{hurault2022proximal}.
For a reference, we also test an unconstrained neural network (baseline NN).
%
%
The network structure given in Definition \ref{def_T}
is applied to the denoisers under study in the current simulation.
The training time is measured by averaging the running time over 30 epochs.

Table \ref{tab:training_speed} gives the training time for images
of $256\times256$ pixels with $8$-layer neural networks.
Owing to the presence of constraints,
all denoisers require longer training time
compared to the unconstrained baseline NN.
The proposed denoiser, however, achieves remarkably
shorter training time compared to MMO and Prox-GS.
Specifically, the proposed denoiser achieves
twenty times faster speed for training
compared to Prox-GS
(twice faster compared to MMO).
We repeat here that the efficiency of the proposed denoiser
in training time comes from the absence of power iteration,
as evidenced now by the training time
for the loss and backward steps shown
in Table \ref{tab:training_speed}.\footnote{
The considerable difference in training time between
MMO and Prox-GS is due to the different operations
in the power iteration.
Specifically, MMO and Prox-GS compute $J_{2T - \Id}(\bm{x})$
and $J_{\nabla \|(\Id - T)(\cdot)\|^2}(\bm{x})$, respectively.
}

Fig.~\ref{fig:training_time} concerns two regimes:
(a) the number of layers in the model changes with the number of pixels
fixed, and
(b) the number of pixels changes with the number of layers fixed.
For convenience, the training time relative to the baseline NN is plotted.
It is observed that the proposed method
consistently shortens the training time significantly across
various configurations.
%

%
\subsection{Deblurring performance}
{\tabcolsep = 3pt
\begin{table*}[t]
  \caption{ Deblurring results in average PSNR [dB] for 30 images from BSDS Dataset}
  \centering
  {\fontsize{9pt}{11pt}\selectfont
  \begin{tabular}{c|ccccccccccc}
    \hline
    \noalign{\vskip 2pt}
    \multirow{2}{*}{\rotatebox{90}{noise level\hspace{10pt}}} & \multirow{2}{*}{\vspace{-43pt}Method} & \multicolumn{8}{c}{Blur kernel}&\multirow{2}{*}{\vspace{-60pt}\it{ Ave.}}\\
    \cline{3-10}
    && \rule{0pt}{2ex} (A) & (B) & (C) & (D) & (E) & (F) & (G) & (H) &\\
    &&
    \includegraphics[width=0.07\linewidth]{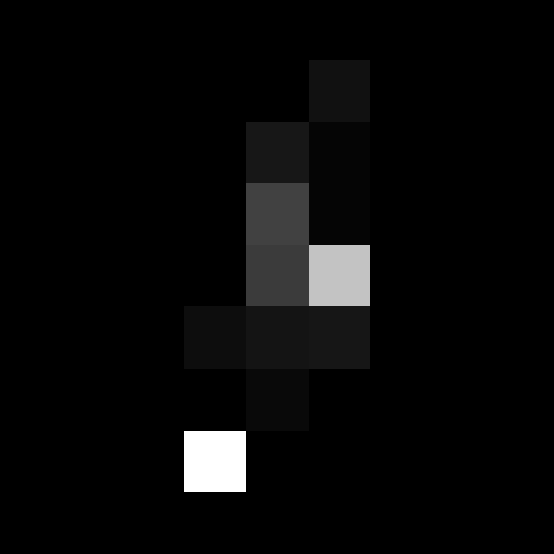}&
    \includegraphics[width=0.07\linewidth]{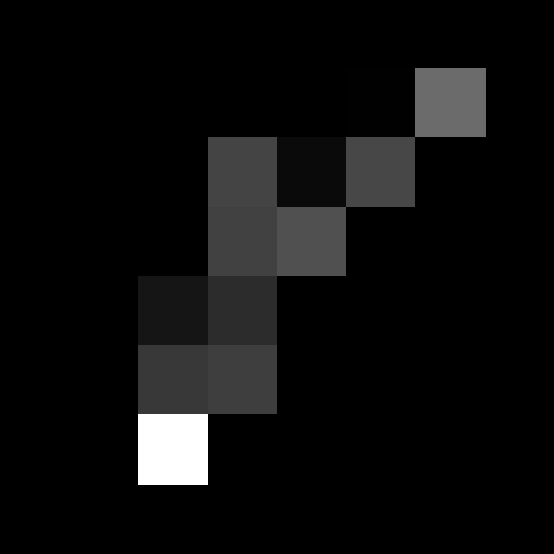}&
    \includegraphics[width=0.07\linewidth]{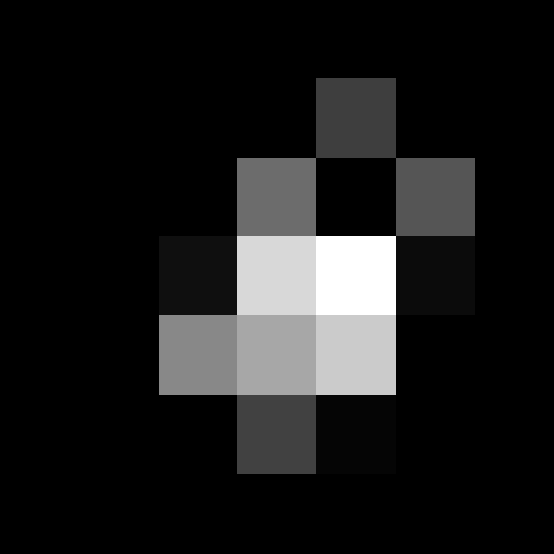}&
    \includegraphics[width=0.07\linewidth]{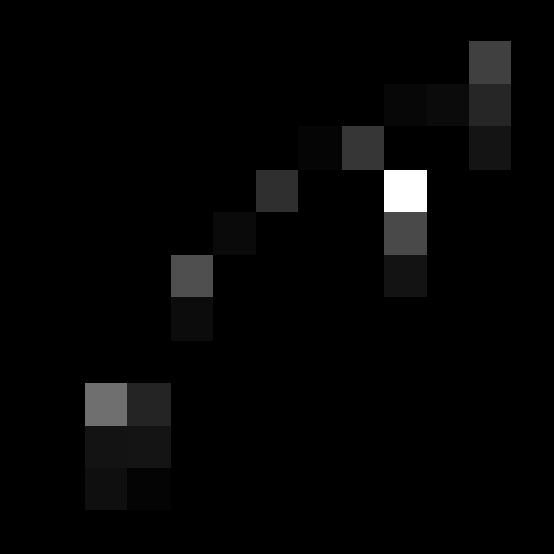}&
    \includegraphics[width=0.07\linewidth]{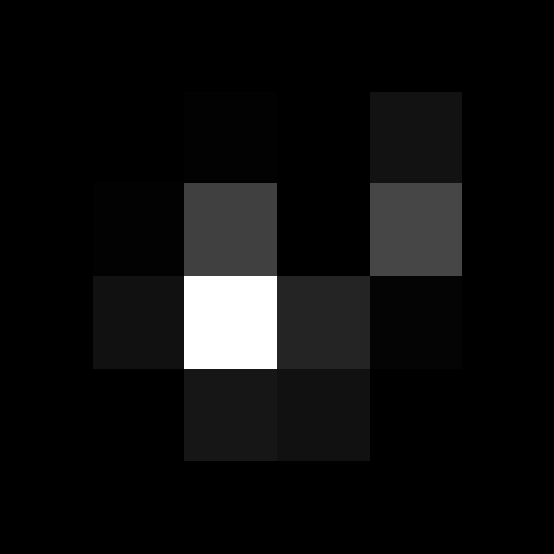}&
    \includegraphics[width=0.07\linewidth]{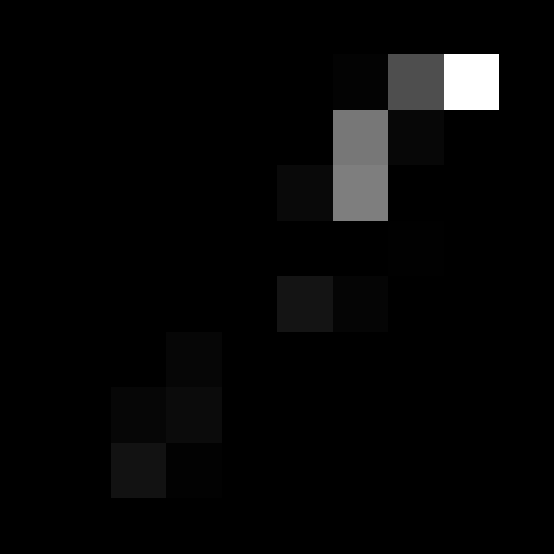}&
    \includegraphics[width=0.07\linewidth]{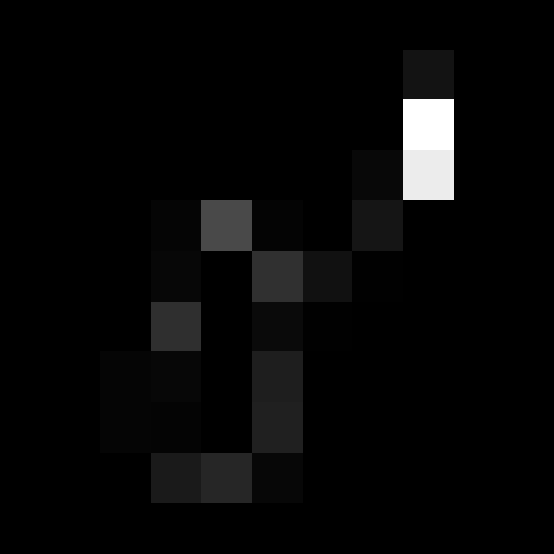}&
    \includegraphics[width=0.07\linewidth]{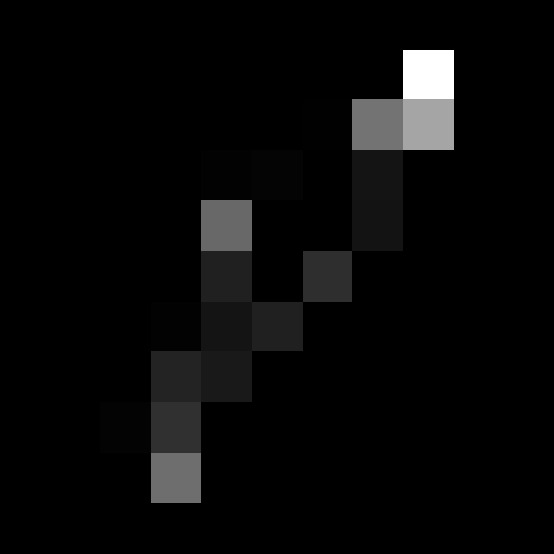}&
    \\
    \hline\hline
    \rule{0pt}{2.5ex}
    \multirow{5}{*}{\rotatebox{90}{$0.01\hspace{2pt}$}}& TV & $32.14$ & $31.31$ & $30.27$ & $31.16$ & $35.24$ & $32.73$ & $32.08$ & $31.27$& $\textit{32.02}$\\
    & MMO & $35.46$ & $34.82$ & $33.54$ & $34.61$ & $38.04$ & $36.05$ & $\bm{35.44}$ & $\bm{34.72}$ & $\textit{35.34}$ \\
    &Prox-GS (DRS) & $\bm{35.68}$ & $\bm{34.93}$ & $\bm{33.74}$ & $\bm{34.75}$ & $\uline{38.05}$ & $\uline{36.18}$ & $35.28$ & $\bm{34.73}$ & $\uline{\textit{35.42}}$ \\[2pt]
    \cline{2-11}
    \rule{0pt}{2ex}
    &MoL-Grad ($\sigma=0.50$) & $35.24$ & $34.65$ & $33.72$ & $\uline{34.64}$ & $\bm{38.48}$ & $36.00$ & $35.27$ & $34.60$ & $\textit{35.33}$\\
    &MoL-Grad ($\sigma=1.75$) & $\uline{35.63}$ & $\uline{34.90}$ & $\bm{33.79}$ & $\bm{34.74}$ & $\bm{38.48}$ & $\bm{36.25}$ & $\uline{35.34}$ & $\uline{34.66}$ & $\textbf{\textit{35.47}}$ \\
    \hline
    \hline
    \rule{0pt}{2.5ex}
    \multirow{5}{*}{\rotatebox{90}{$0.03\hspace{2pt}$}}& TV & $26.92$ & $26.20$ & $26.53$ & $26.15$ & $28.98$ & $27.33$ & $27.13$ & $26.49$& $\textit{26.97}$\\
    & MMO & $30.29$ & $\uline{29.53}$ & $29.60$ & $29.33$ & $\uline{32.38}$ & $30.67$ & $30.37$ & $29.74$ & $\textit{30.24}$\\
    &Prox-GS (DRS) & $\bm{30.53}$ & $\bm{29.86}$ & $\uline{29.78}$ & $\bm{29.77}$ & $32.18$ & $\uline{30.86}$ & $\uline{30.65}$ & $\bm{30.12}$ & $\uline{\textit{30.47}}$ \\[2pt]
   \cline{2-11}
   \rule{0pt}{2ex}
    &MoL-Grad ($\sigma=0.50$) & $29.87$ & $29.47$ & $29.66$ & $29.43$ & $\bm{32.70}$ & $30.84$ & $30.58$ & $29.92$ & $\textit{30.31}$\\
&MoL-Grad ($\sigma=1.75$) & $\uline{30.43}$ & $\bm{29.87}$ & $\bm{29.83}$ & $\uline{29.73}$ & $\bm{32.70}$ & $\bm{31.09}$ & $\bm{30.69}$ & $\uline{30.09}$ & $\textbf{\textit{30.55}}$ \\
    \hline
  \end{tabular}
  }
  \label{tab:psnr_deblur}
\end{table*}
}

%

The performance of the proposed denoiser is
compared to that of total variation (TV),
MMO\cite{pesquet2021_MMO}, and Prox-GS\cite{hurault2022proximal}.
To assess the deblurring performance,
30 images from the Berkeley Segmentation Dataset (BSDS300) \cite{BSDS}
are used with the real-world motion kernels from
\cite{Levin_blurkernel}, scaled down by a factor of 1/2.
We set $\bm{v}_0=\bm{0}$ and $\bm{u}_0=\bm{0}$ in Algorithm \ref{alg:primal_dual}.
The algorithm terminates when
$\|\bm{v}_k-\bm{v}_{k-1}\|^2/\|\bm{v}_0\|^2 < 10^{-10}$ or if
the iteration number exceeds 500.

%
%
%
In our experiment, the estimated Lipschitz constant was $\hat{L}_D=2.28$
(see Remark \ref{remark:LD}), to which the corresponding $\sigma$ is
$\sigma=0.78$.
%
We consider the two settings as mentioned in Remark \ref{remark:LD}:
(i) $\sigma=0.50$ (the conservative setting by keeping in mind that
$L_D\geq \hat{L}_D$), and
(ii) $\sigma = 1.75$ (the relaxed setting to assign a larger weight to $\varphi$).
For MMO and Prox-GS, publicly available pretrained models are used.
In addition, Prox-GS is evaluated using only the
Douglas-Rachford splitting (DRS) algorithm.

\begin{figure}[t]
    \centering
    \includegraphics[width=0.9\linewidth]{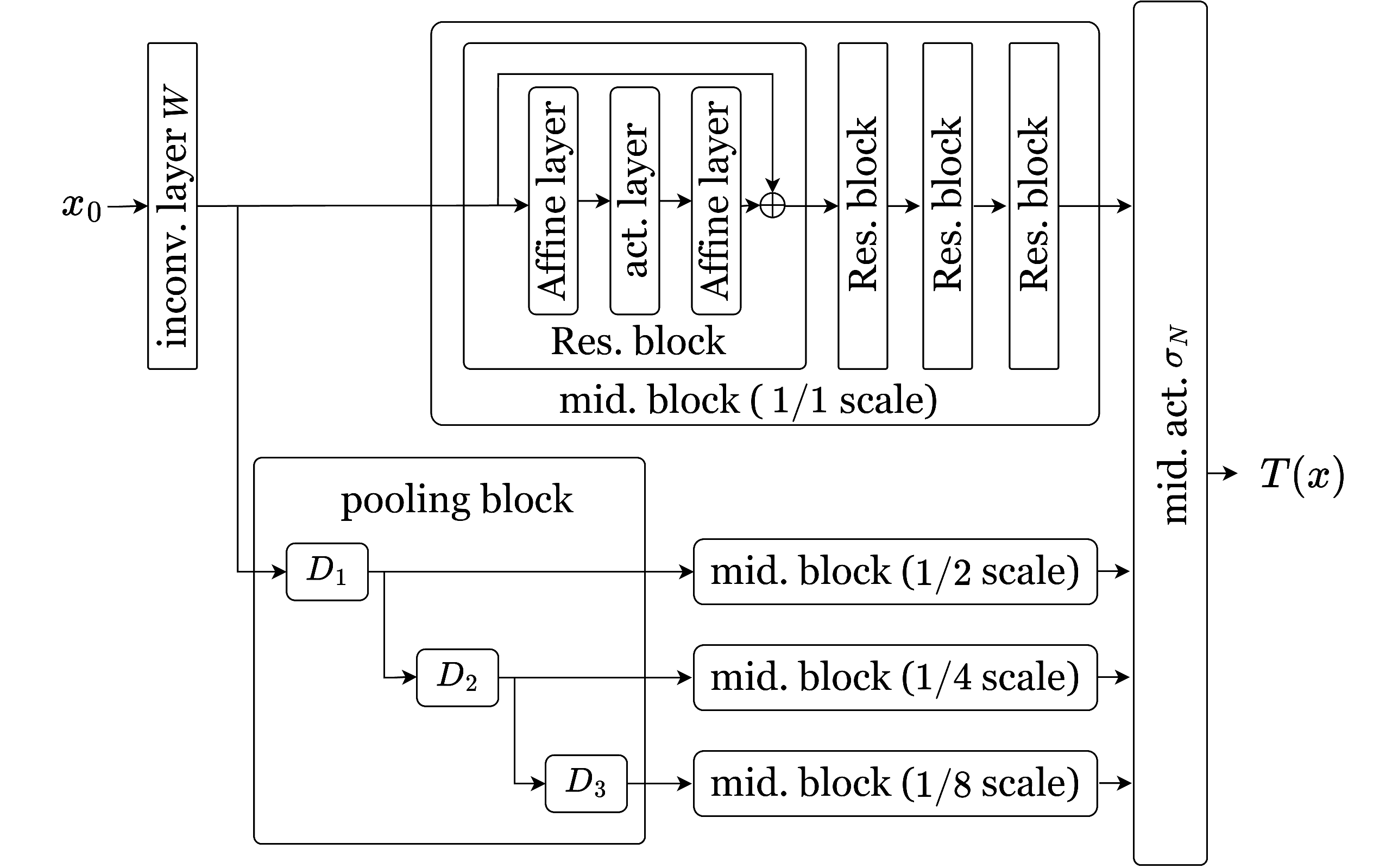}
    \caption{The architecture of $T$ used in the deblurring results.
    An affine layer represents $\bm{W} \cdot + \bm{b}$, where $\bm{W}$ is a convolutional layer.
    Each $D_1$, $D_2$, and $D_3$ represents $\bm{D}\circ\bm{W}$, where $\bm{D}$ is $2 \times 2$ average pooling.}
    \label{fig:structure_T}
\end{figure}

\begin{figure}[tbp]
  \centering
  \subfloat[True image]{\includegraphics[width=0.49\linewidth]{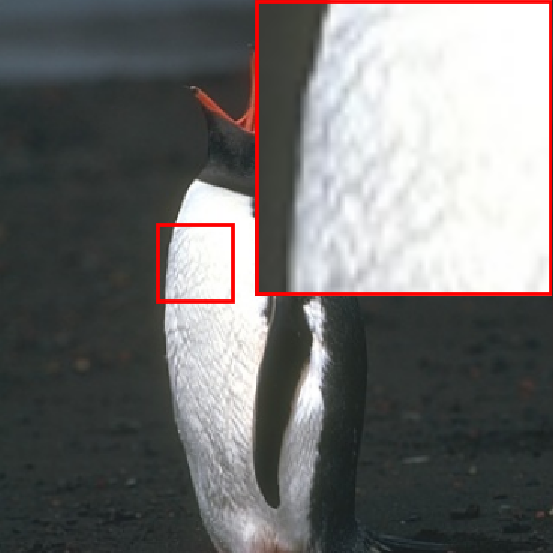}}
  \hfill
  \subfloat[Blurred image (PSNR $=26.62$)]{\includegraphics[width=0.49\linewidth]{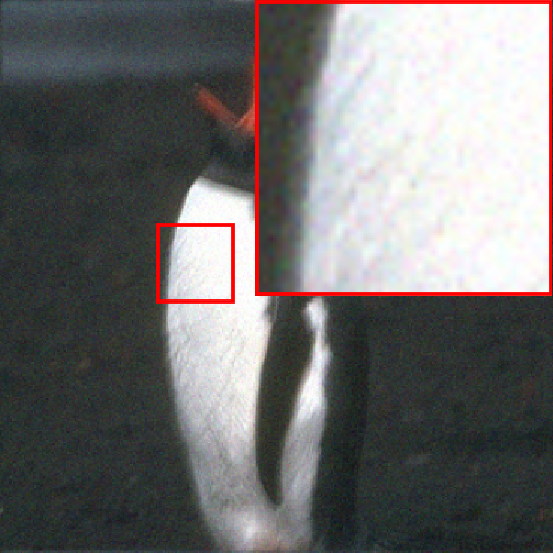}}

  \subfloat[TV (PSNR $=32.12$)]{\includegraphics[width=0.49\linewidth]{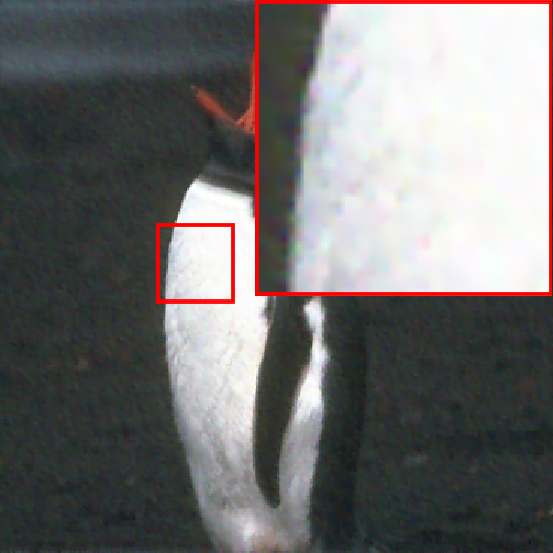}}
  \hfill
  \subfloat[MMO (PSNR $=35.04$)]{\includegraphics[width=0.49\linewidth]{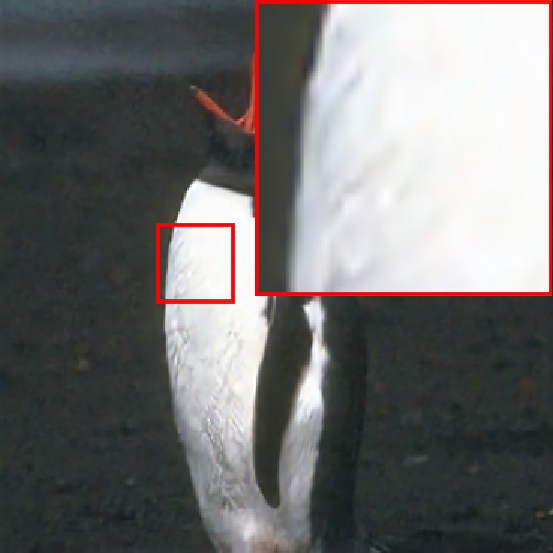}}

  \subfloat[Prox-GS (PSNR $=33.92$)]{\includegraphics[width=0.49\linewidth]{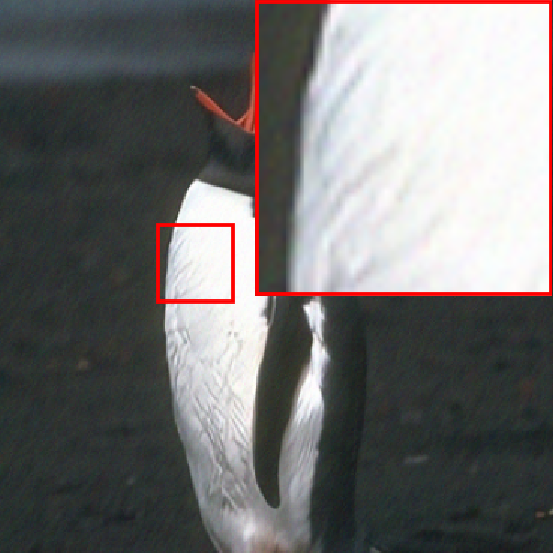}}
  \hfill
  \subfloat[MoL-Grad (Proposed)\\ \hspace{10pt}(PSNR $=35.32$)]{\includegraphics[width=0.49\linewidth]{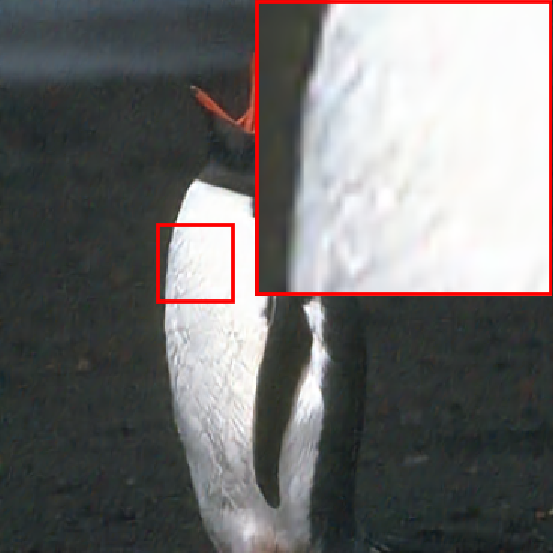}}
  \caption{Deblurring results of the penguin image (blur kernel (B), noise level $\sigma_{\varepsilon}=0.03$) (PSNR in dB).}
  \label{fig:deblurring_penguin}
\end{figure}

\begin{figure}[tbp]
  \centering
  \subfloat[True image]{\includegraphics[width=0.49\linewidth]{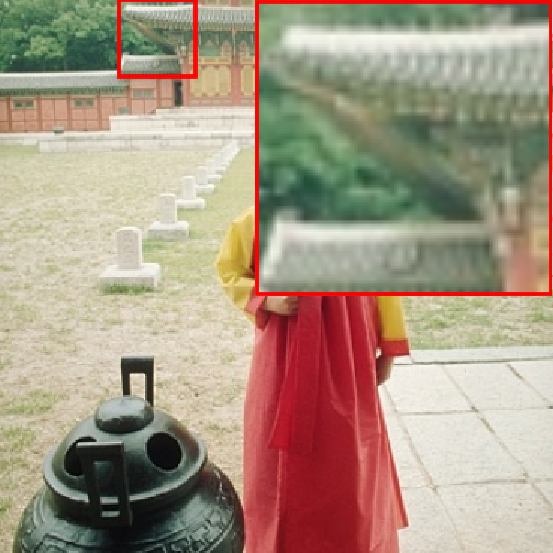}}
  \hfill
  \subfloat[Blurred image (PSNR $=19.28$)]{\includegraphics[width=0.49\linewidth]{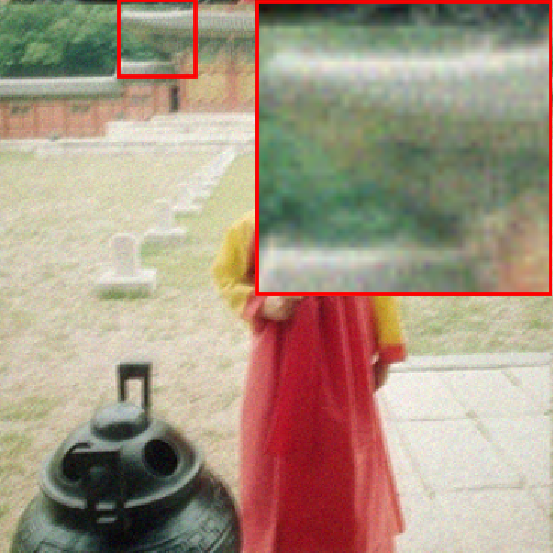}}

  \subfloat[TV (PSNR $=28.18$)]{\includegraphics[width=0.49\linewidth]{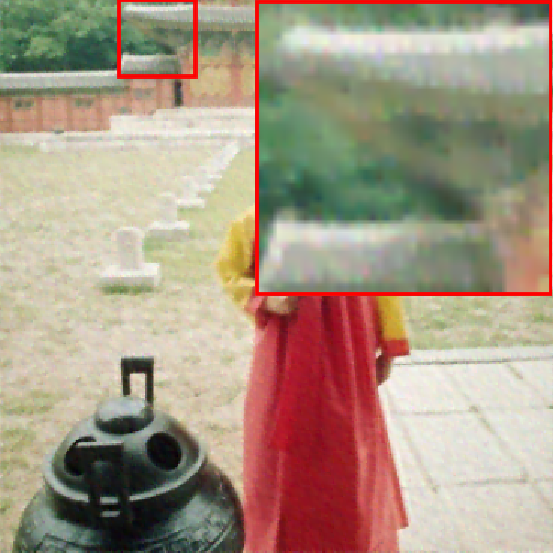}}
  \hfill
  \subfloat[MMO (PSNR $=30.70$)]{\includegraphics[width=0.49\linewidth]{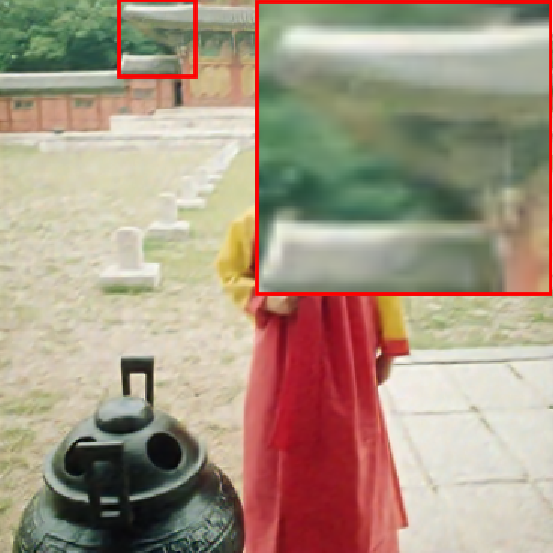}}

  \subfloat[Prox-GS (PSNR $=31.24$)]{\includegraphics[width=0.49\linewidth]{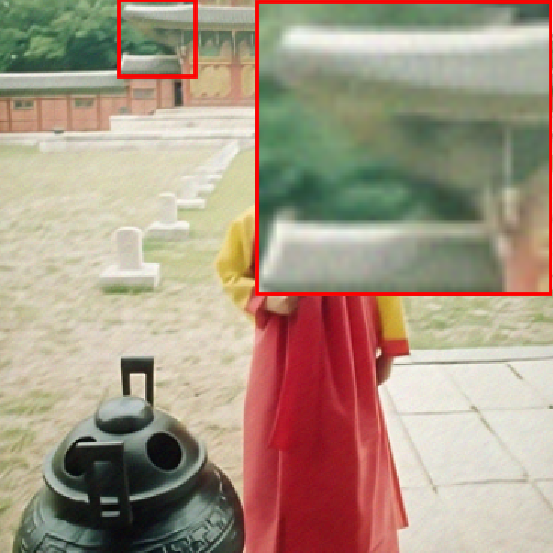}}
  \hfill
  \subfloat[MoL-Grad (Proposed)\\ \hspace{10pt} (PSNR $=31.52$)]{\includegraphics[width=0.49\linewidth]{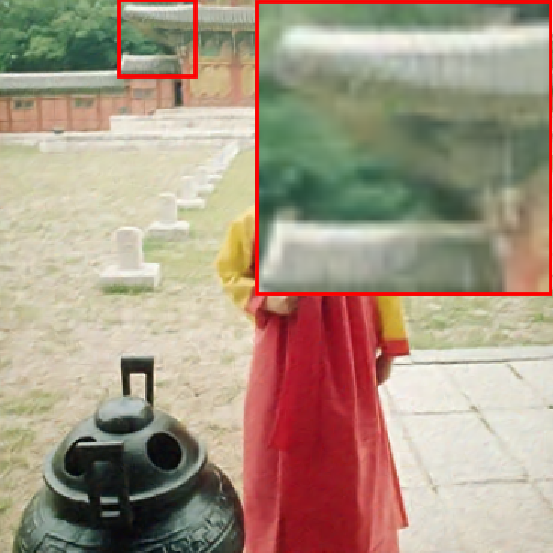}}
  \caption{Deblurring results of the palace image (blur kernel (F), noise level $\sigma_\varepsilon=0.03$) (PSNR in dB).}
  \label{fig:deblurring_palace}
\end{figure}

Residual blocks and pooling blocks depicted in
Fig.~\ref{fig:structure_T}
are considered in this experiment.
We set $\sigma_{n, i}$ to a smoothed ReLU (sReLU) defined by
\begin{equation}
\text{sReLU}(x) :=
    \begin{cases}
    x &\text{if $\gamma < x$}, \\
    x^2/(4\gamma) + x/2 + \gamma/4 &\text{if $-\gamma \leq x \leq \gamma$}, \\
    0 &\text{otherwise},
    \end{cases}
\end{equation}
for $\gamma\in\mathbb{R}_{++}$.
This setting of $\sigma_{n,i}$ does not violate
assumptions \ref{assumption_conv_sigma} and \ref{assumption_lip_sigma}
and condition \ref{assumption_lip_sigmaprime} of Theorem \ref{theorem_mol_NN}.
For training, images are randomly cropped to $128 \times 128$ from
Flickr2K \cite{Flicker2K}, DIV2K \cite{Agustsson2017_DIV2K}, and
Waterloo Exploration Database \cite{ma2017_waterloo}
by following the way of \cite{Zang2022_drunet},
\cite{hurault2022_gradient};
the noise level is set to $\sigma_\varepsilon = 0.15$.
In addition, the loss function given in \eqref{eq:def_loss} is optimized
using Adam~\cite{Kingma2014_Adam}.
The learning rate is set to $1 \times 10^{-4}$ for the first 2700
training epochs and
to $2.5\times10^{-6}$ for the remaining 300 epochs.
It takes approximately two weeks to train the model on the NVIDIA A30 GPU.
Table~\ref{tab:psnr_deblur} summarizes the deblurring results.
The proposed method preserves performance comparable to that of
those state-of-the-art methods.
Figures~\ref{fig:deblurring_penguin} and \ref{fig:deblurring_palace}
display specific instances of the deblurred images.
It can be seen that the proposed method restores the fine patterns
better, such as the texture on the penguin's belly and
the stripe pattern on the roof.
From these results,
the proposed method successfully achieves the state-of-the-art performance
in the practical application.
We stress that this excellent performance comes with
the significant reduction of training time (see Section \ref{subsec:training_speed}).

\section{Conclusion}
\label{sec:conclusion}
In this paper,
we presented the weight-tying multi-layer nonnegative neural network architecture
yielding a MoL-Grad denoiser that is a key component
of the proposed image restoration algorithm.
The architecture was derived from the gradient of the sum of
the smooth convex functions each of which
is associated with the composition of the
activation function and an affine operator.
The weights are therefore tied between the encoder and decoder,
and the monotonicity of the denoiser
is ensured owing to the weight nonnegativity.
The results were extended to the architecture including skip
connections which are widely used in image restoration tasks.
The MoL-Grad denoiser based on the proposed architecture
is free from Lipschitz constraints
in sharp contrast to the previous theoretical studies on
convergence of the PnP method.
This leads to the significant advantages over the nonexpansive denoisers
in terms of the training time as well as the denoising performance,
as demonstrated by the simulations.
In particular, the proposed method achieved the state-of-the-art deblurring performance among those PnP algorithms having convergence guarantees.
As the ``implicit'' regularizer induced by the proposed denoiser tends
to be weakly convex,
the additional penalty of the projection-involving quadratic function
was introduced.
Because of that, even in the ill-conditioned case,
the overall convexity is preserved,
and the sequence of vectors generated by the primal-dual splitting algorithm
converges to a minimizer of the cost function involving the implicit regularizer.
The present study showed the first practical multi-layer neural network architecture
yielding a MoL-Grad denoiser.
We emphasize the key point once again that
the proposed neural network empowered approach is ``explainable''
in the sense of disclosing the objective minimized.
This work suggests that the PnP method adopting MoL-Grad denoisers
will be useful in real-world applications.
There still remains plenty of room, nevertheless,
for finding better performing sub-families of the MoL-Grad denoiser.
\appendices
\renewcommand{\theequation}{A.\arabic{equation}}
\setcounter{equation}{0}
\section{proof of lemma~\ref{claim_conv}}
\label{app:proof_pf_lemma1}
Let $\alpha \in (0, 1)$.
By assumption, since $t_i$ is convex and $w_i\in\mathbb{R}_+$, we have
\begin{equation}
    w_i t_i (\alpha \bm{x} + (1-\alpha) \bm{y}) \leq w_i \left[ \alpha t_i (\bm{x})+(1-\alpha)t_i(\bm{y}) \right],
\end{equation}
for every $\bm{x}, \bm{y} \in \mathbb{R}^n$.
In addition, we have
\begin{align}
    &\bm{w}^\mathsf{T} T(\alpha \bm{x} + (1-\alpha) \bm{y}) =\sum_{i=1}^q w_{i} t_i(\alpha \bm{x} +(1-\alpha) \bm{y}) \notag\\
    &\hspace{30pt}\leq \sum_{i=1}^q w_{i} \left[ \alpha t_i(\bm{x})+ (1-\alpha)t_i(\bm{y})\right]\notag\\
    &\hspace{30pt}=\alpha \sum_{i=1}^q \left[ w_{i}t_i(\bm{x})\right] + (1-\alpha) \sum_{i=1}^q \left[ w_{i} t_i(\bm{y})\right]\notag\\
    &\hspace{30pt}= \alpha \bm{w}^\mathsf{T} T(\bm{x}) + (1-\alpha) \bm{w}^\mathsf{T} T(\bm{y}),\label{eq:lemma1_1}
\end{align}
which implies that the function $\bm{x}\mapsto\bm{w}^\mathsf{T}T(\bm{x})$ is convex.
Since $\sigma$ is non-decreasing and convex by assumption, it follows from \eqref{eq:lemma1_1} that
\begin{align}
    &\sigma(\bm{w}^\mathsf{T} T(\alpha \bm{x} + (1-\alpha) \bm{y})+b) \notag\\
    &\hspace{15pt} \leq \sigma(\alpha \bm{w}^\mathsf{T} T(\bm{x}) + (1-\alpha) \bm{w}^\mathsf{T} T(\bm{y}) +b) \notag\\
    &\hspace{15pt} \leq \alpha \sigma(\bm{w}^\mathsf{T} T(\bm{x})+b) + (1-\alpha) \sigma(\bm{w}^\mathsf{T} T(\bm{y})+b),
\end{align}
which implies $\sigma(\bm{w}^\mathsf{T} T(\cdot) +b)$ is convex.
\qed
\renewcommand{\theequation}{B.\arabic{equation}}
\setcounter{equation}{0}
\section{proof of lemma~\ref{memo_lip2}}
\label{app:proof_pf_lemma2}
For every $\bm{x}, \bm{y} \in \mathbb{R}^s$, we have
\begin{align}
    &\big\|F(\bm{x})-F(\bm{y})\big\| \notag\\
    & =\big\| \bm{W}\big(T(\bm{x})\odot R(\bm{x})\big) - \bm{W}\big(T(\bm{y})\odot R(\bm{y})\big)\big\| \notag\\
    &\leq \|\bm{W}\| \| \big(T(\bm{x})\odot R(\bm{x})\big) - \big(T(\bm{y})\odot R(\bm{y})\big)\| \notag \\
    &=\|\bm{W}\| \big\| \big( T(\bm{x}) \odot R(\bm{x})\big) - \big(T(\bm{y}) \odot R(\bm{x})\big) \notag\\
    & \hspace{55pt}+ \big(T(\bm{y})\odot R(\bm{x})\big) - \big(T(\bm{y})\odot R(\bm{y})\big)\big\| \notag\\
    & \leq \|\bm{W}\|  \Big( \big\|\big(T(\bm{x})\odot R(\bm{x})\big) - \big(T(\bm{y})\odot R(\bm{x})\big) \big\| \notag\\
    &\hspace{55pt}+ \big\| \big(T(\bm{y})\odot R(\bm{x})\big)- \big(T(\bm{y})\odot  R(\bm{y})\big)\big\| \Big) \notag\\
    %
    &= \|\bm{W}\| \Big( \big\|\big(T(\bm{x})-T(\bm{y}) \big)\odot  R(\bm{x}) \big\| \notag\\
    &\hspace{55pt}+ \big\|T(\bm{y}) \odot \big(R(\bm{x})-R(\bm{y})\big)\big\|\Big)\notag\\
    &\leq \|\bm{W}\| \Big(\max \big(  \big| R(\bm{x})\big| \big) \big\|T(\bm{x})-T(\bm{y})\big\| \notag\\
    &\hspace{55pt}+ \max \big( \big| T(\bm{y})\big| \big) \big\|R(\bm{x})-R(\bm{y}) \big\| \Big)\notag\\
    &\leq \|\bm{W}\| \Big( \max \big( \big| R(\bm{x})\big| \big) L_T \notag\\
    &\hspace{55pt} +\max \big( \big| T(\bm{y}) \big| \big) L_{ R} \Big) \| \bm{x}-\bm{y}\| \notag\\
    &=L_F \|\bm{x}-\bm{y}\|,
\end{align}
where $\max(\lvert \bm{x} \rvert)$  denotes the maximum absolute value among the components of $\bm{x}$, and $L_T, L_{R}$ are Lipschitz constants of $T$ and $R$, respectively.
We obtain the third inequality by applying the fact that
$\|\bm{x}\odot\bm{y} \| \leq \| \max\left(|\bm{x}|\right) \bm{1}\odot  \bm{y} \|=\max (|\bm{x}|) \|\bm{y}\|$, where $\bm{1}:=[1,  1, \cdots, 1]^\mathsf{T}\in\mathbb{R}^q$.
Hence, $F$ is $L_F$-Lipschitz continuous with a constant $L_F \leq \|\bm{W}\|(M_r L_T+M_t L_R)$.
Finally, the boundedness of $f_i$ can be verified by $|f_i(\bm{x})| =|[\bm{W}]_i \big(T(\bm{x})\odot R(\bm{x})\big)|\leq\sum_{j=1}^q |[\bm{W}]_{i, j} |M_t M_r$ for every $\bm{x}\in\mathbb{R}^s$. \qed
\renewcommand{\theequation}{C.\arabic{equation}}
\setcounter{equation}{0}
\section{proof of theorem \ref{theorem_mol_NN}}
\label{app:proof_of_theorem1}
\begin{figure*}[t]
    \centering
    \includegraphics[width=1\linewidth]{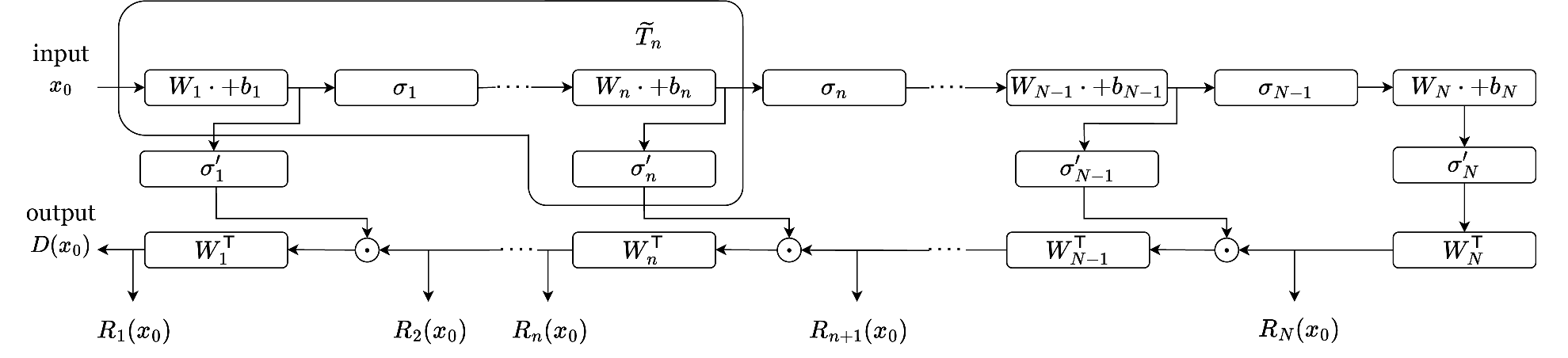}
    \caption{Definitions of $\widetilde{T}_n$ and $R_n$.}
    \label{fig:model_of_R}
\end{figure*}
{\normalfont (a)\hspace{5pt}}
Define the function $t_{n:1, i} : \bm{x} \mapsto (T_{n:1} (\bm{x}))_i$ for $n>1$, where $(T_{n:1}(\bm{x}))_i$ is $i$th component of the vector $T_{n:1}(\bm{x})$.
Likewise, let $b_{n, i}$ be the $i$th component of $\bm{b}_n$.

The proof is based on mathematical induction.
For the first layer, for every $i\in\{1, 2, \cdots, d_{1}\}$, the function
\begin{equation}
  t_{1, i}: \bm{x}\mapsto\sigma_{1, i}([\bm{W}_1]_i \bm{x}+ b_{1, i})
\end{equation}
is convex because the function $\sigma_{1, i}$ is convex by assumption \ref{assumption_conv_sigma}.
Assume that $t_{n:1, j}$ is convex for every $j\in\{1, 2, \cdots, d_n\}$.
Then, for the $(n+1)$st layer, we have
\begin{align}
  t_{n+1:1, i} (\bm{x}) &= \sigma_{n+1, i}([\bm{W}_{n+1}]_i T_{n:1}(\bm{x}) + b_{n+1, i})\notag\\
  &= \sigma_{n+1, i}\left[\sum_{j=1}^{d_n} [\bm{W}_{n+1}]_{i, j} t_{n:1, j}(\bm{x}) + b_{n+1, i}\right]
\end{align}
for every $i\in\{1, 2, \cdots, d_{n+1}\}$.
Thus, by Lemma~\ref{claim_conv} together with assumptions \ref{assumption_conv_sigma} and \ref{assumption_conv_W}, $t_{n+1:1, i}$ is also convex.
By induction, it can be verified that $t_i = t_{N:1, i}$ is convex (and clearly differentiable as well). \\
{\normalfont (b)\hspace{5pt}}
Define (see Fig.~\ref{fig:model_of_R})
\begin{align}
\widetilde{T}_n(\bm{x}) :\hspace{-2pt}&= \left[ \tilde{t}_{n,1}(\bm{x}), \tilde{t}_{n,2}(\bm{x}), \cdots, \tilde{t}_{n,d_n}(\bm{x})\right]^\mathsf{T}\\
:\hspace{-2pt}&=\bm{\sigma}_n'(\bm{W}_n T_{n-1:1}(\bm{x})+\bm{b}_n).\label{eq:def_tilde_T}
\end{align}
According to assumption \ref{assumption_lip_sigmaprime}, the Lipschitz continuity of $\sigma'_{n,i}$ for every $i$ guarantees that the operator $\sigma'_n$ is Lipschitz continuous.
Furthermore, $\bm{\sigma}_n$ is Lipschitz continuous by assumption \ref{assumption_lip_sigma}, and hence its composition $T_{n-1:1}$ with affine mappings.
Hence, $\widetilde{T}_n$ defined by \eqref{eq:def_tilde_T} is Lipschitz continuous.
Moreover, $\tilde{t}_{n,i}$ is also bounded, as assumption \ref{assumption_lip_sigma} guarantees the boundedness of $\sigma'_{n,i}$.
To facilitate the proof, we additionally define the operator $R_n$ recursively from $n = N$ to $n = 1$, and then define the functions $r_{n,i}$ as follows:
\begin{align}
    R_n (\bm{x}) &:=\left[ r_{n, 1}(\bm{x}), r_{n, 2}(\bm{x}), \cdots, r_{n, d_n}(\bm{x})\right]^\mathsf{T}\\
    &:=
    \begin{cases}
        \bm{W}^\mathsf{T}_N \widetilde{T}_N(\bm{x}) &\text{if $n=N$}\\
        \bm{W}^\mathsf{T}_{n}\big[ \widetilde{T}_{n}(\bm{x}) \odot \left( R_{n+1}(\bm{x}) \right)\big] &\text{if $n< N$.}
    \end{cases}
\end{align}
Using this notation, the expression for $D(\bm{x}_0)$ given in \eqref{eq:express_D} can be rewritten in terms of $R_n$.
First, $D$ can be expressed as
\begin{align}
    D(\bm{x}_0) &=\left[ \prod_{n=1}^{N-1}  \mathrm{J}_{T_n} (x_{n-1})^\mathsf{T} \right]\notag\\
    &\hspace{30pt}\circ \bm{W}_N^\mathsf{T}\circ \bm{\sigma}'_N \circ (\bm{W}_N \cdot+\bm{b}_N)\circ T_{N-1:1} (\bm{x}_0)\notag\\
    &=\left[ \prod_{n=1}^{N-1}  \mathrm{J}_{T_n} (x_{n-1})^\mathsf{T} \right]\circ \bm{W}_N^\mathsf{T}\widetilde{T}_N (\bm{x}_0)\notag\\
    & = \left[ \prod_{n=1}^{N-1}  \mathrm{J}_{T_n} (x_{n-1})^\mathsf{T} \right]\circ R_N(\bm{x}_0).\notag\\
    &=\left[ \prod_{n=1}^{N-2}  \mathrm{J}_{T_n}(x_{n-1})^\mathsf{T}\right]\circ
    \left[\mathrm{J}_{T_{N-1}}(\bm{x}_{N-2})^\mathsf{T} R_N(\bm{x}_0)\right].\label{eq:express_D_using_R}
\end{align}
Here, using \eqref{eq:express_jaco}, we have
\begin{align}
    &\mathrm{J}_{T_{N-1}}(\bm{x}_{N-2})^\mathsf{T}R_N(\bm{x}_0)\notag\\
    &= \bm{W}^\mathsf{T}_{N-1} \big(\bm{\sigma}'_{N-1} (\bm{W}_{N-1}\bm{x}_{N-2}+\bm{b}_{N-1})\odot R_N(\bm{x}_0)\big) \notag\\
    &=\bm{W}^\mathsf{T}_{N-1} \big(\widetilde{T}_{N-1} (\bm{x}_0) \odot R_N(\bm{x}_0)\big)\notag\\
    &=R_{N-1}(\bm{x}_0), \label{eq:leading_RN-1}
\end{align}
where the second equality can be verified by \eqref{eq:def_x_n} and \eqref{eq:def_tilde_T}.
Similarly, for every $n$, we have
\begin{align}
    \mathrm{J}_{T_{n-1}}(\bm{x}_{n-2})^\mathsf{T}R_n(\bm{x}_0)=R_{n-1}(\bm{x}_0).\label{eq:leading_Rn}
\end{align}
Finally, combining \eqref{eq:express_D_using_R} and \eqref{eq:leading_RN-1} and then applying \eqref{eq:leading_Rn}, we obtain (see Fig.~\ref{fig:model_of_R})
\begin{align}
    D(\bm{x}_0)&=\left[ \prod_{n=1}^{N-2}  \mathrm{J}_{T_n}(\bm{x}_{n-1})^\mathsf{T}\right]\circ
     R_{N-1}(\bm{x}_0)\notag\\
     &=\left[ \prod_{n=1}^{N-3}  \mathrm{J}_{T_n}(\bm{x}_{n-1})^\mathsf{T}\right]\circ R_{N-2}(\bm{x}_0)\notag\\
     &\hspace{25pt}\vdots\notag\\
     &=\mathrm{J}_{T_1}(\bm{x}_{0})^\mathsf{T}R_2(\bm{x}_0)\notag\\
     &=R_1(\bm{x}_0).
\end{align}
To prove the Lipschitz continuity of $D$, we show by induction that each $R_n$ is Lipschitz continuous and each function $r_{n,i}$ is bounded, using the Lipschitz continuity of $\widetilde{T}_n$ and the boundedness of $\tilde{t}_{n,i}$.
First of all, $R_N=\bm{W}_N^\mathsf{T}\widetilde{T}_N$ is Lipschitz continuous and $r_{N, i}=[\bm{W}_N^\mathsf{T}]_i\tilde{T}_{N}=\sum_{j=1}^{d_n}[\bm{W}_N^\mathsf{T}]_{i, j}\tilde{t}_{N, j}$ is bounded.
Suppose now, for $n<N$, that $R_{n+1}$ is Lipschitz continuous and that $r_{n+1, i}$ is bounded.
Then, since $\widetilde{T}_n$ and $\tilde{t}_{n, i}$ are also Lipschitz continuous and bounded, $R_{n}$ and $r_{n, i}$ are Lipschitz continuous and bounded, respectively, by Lemma~\ref{memo_lip2} under the correspondence $F:=R_n$, $R:=R_{n+1}$, $T:=\widetilde{T}_n$, and $\bm{W}:=\bm{W}_{n}^\mathsf{T}$.
It therefore follows that $R_1=D$ is Lipschitz continuous.
\\
{\normalfont (c)\hspace{5pt}}
By (a), assumptions \ref{assumption_conv_sigma} and \ref{assumption_conv_W} imply that the function $\psi:=\sum_{i=1}^{d_N} t_i$ is convex and differentiable.
In addition, by (b), assumptions \ref{assumption_lip_sigma} and \ref{assumption_lip_sigmaprime} imply that the denoiser $D=\nabla\psi$ is Lipschitz continuous.
The assumption $L_D > 1$ implies $\beta \in (0, 1)$.
Hence, Fact~\ref{fact:def_molgrad} verifies that $D$ is the MoL-Grad denoiser.
\qed
\renewcommand{\theequation}{D.\arabic{equation}}
\setcounter{equation}{0}
\section{proof of theorem \ref{theorem_skip}}
\label{app:proof_of_theorem2}
\noindent{\normalfont (a)\hspace{5pt}} We observe that
    \begin{align}
        T^{b\leftarrow a} &= T_{N:b} \circ (T_{a:1}+T_{b-1:1}).
    \end{align}
    Under assumptions \ref{assumption_conv_sigma} and \ref{assumption_conv_W} for $n\in\{1, 2, \cdots, a, \cdots, b-1\}$, the functions $t_{a-1:1, i}$ and $t_{b:1, i}$ are convex, and hence $t_{a-1:1, i} + t_{b:1, i}$ is also convex.
    Under this observation, together with assumptions \ref{assumption_conv_sigma} and \ref{assumption_conv_W} for $n=a$, Lemma~\ref{claim_conv} implies convexity of a function $\bm{x}\mapsto\sigma_{a, i} ([\bm{W}_a]_{i} (T_{a-1:1}(\bm{x})+T_{b:1}(\bm{x})))$.
    Hence, the convexity of $t^{b\leftarrow a}_i$ is obtained in the same way as in the proof of Theorem~\ref{theorem_mol_NN}(a).\\
\noindent{\normalfont (b)\hspace{5pt}}
Define
\begin{align}
    \tsup{T}_n (\bm{x}) &:= \left[ \tsup{t}_{n, 1}(\bm{x}), \tsup{t}_{n, 2} (\bm{x}), \cdots, \tsup{t}_{n, d_n} (\bm{x})\right]^\mathsf{T}\\
    &=
    \begin{cases}
        \widetilde{T}_n (\bm{x}) &\text{if $n < b$}\\
        \bm{\sigma}_n'  \big[ \bm{W}_n (\Id+T_{b-1:a+1}) &\\
        \hspace{20pt}\circ T_{a:1}(\bm{x}) +\bm{b}_n \big] &\text{if $n=b$} \\
        \bm{\sigma}_n' \big[ \bm{W}_n T_{n-1:b} \circ(\Id+T_{b-1:a+1}) & \\
        \hspace{20pt}\circ T_{a:1} (\bm{x}) + \bm{b}_n \big] &\text{if $n> b$.}
    \end{cases}
\end{align}
For $n\in\{0, 1, 2, \cdots N\}$, we define the $n$th layer output
(equivalently, the $(n+1)$st layer input) by
\begin{align}
    \hspace{-4pt}\bm{x}^{b\leftarrow a}_{n} :=
    \begin{cases}
    \bm{x}_0 & \text{if $n=0$}\\
    T_{n:1}(\bm{x}_0)                &\text{if $n < b-1$}\\
    \big(T_{b-1:1}+T_{a:1}\big)(\bm{x}_{0})   &\text{if $n=b-1$}\\
    T_{n:b}\circ\big(T_{b-1:1}+T_{a:1}\big)(\bm{x}_{0}) &\text{if $n> b-1$},
    \end{cases}
\end{align}
for $n\in\{0, 1,2,\cdots N\}$.
Using this, we obtain $\tsup{T}_n(\bm{x}_0^{b\leftarrow a})=\bm{\sigma}'_n(\bm{W}_n \bm{x}_{n-1}^{b\leftarrow a}+\bm{b}_n)$ for every $n=\{ 1, 2, \cdots, N\}$.
We also define $\widetilde{R}_n$ and $\widetilde{R}^{b\leftarrow a}_{ n}$ recursively from $n=N$ to $n=1$ as follows:
\begin{align}
    \widetilde{R}_n (\bm{x}) &:=
    \begin{cases}
        \bm{W}_N^\mathsf{T} \tsup{T}_N (\bm{x}) &\text{if $n=N$}\\
        \bm{W}_n^\mathsf{T} \left[ \tsup{T}_n (\bm{x})\odot \left(\widetilde{R}_{n+1}(\bm{x}) \right) \right] &\text{if $n< N$}
    \end{cases}
\end{align}
\vspace{-15pt}
\begin{align}
    \widetilde{R}^{b\leftarrow a}_{ n}&(\bm{x}) \notag \\
    :\hspace{-2pt}=&
    \begin{cases}
        \bm{W}_N^\mathsf{T}\tsup{T}_N (\bm{x}) &\text{if $n =N$}\\
        \widetilde{R}^{b\leftarrow a}_{ n+1} (\bm{x}) &\text{if $a < n < b$}\\
        \bm{W}^\mathsf{T}_n \left[\tsup{T}_{n} (\bm{x}) \odot \left( \widetilde{R}^{b\leftarrow a}_{ n+1} (\bm{x})\right)\right]& \text{otherwise}.
    \end{cases}
\end{align}
Then, using the same approach as in the proof of Theorem~\ref{theorem_mol_NN}(b), it follows that
\begin{align}
    &D^{b\leftarrow a}(\bm{x}_0^{b\leftarrow a}) \nonumber\\
    &= \mathrm{J}_{T_{a:1}}(\bm{x}^{b\leftarrow a}_{ 0})^\mathsf{T} \left( \bm{I}+\mathrm{J}_{T_{b-1:a+1}}(\bm{x}^{b\leftarrow a}_{a})^\mathsf{T} \right) \mathrm{J}_{T_{N-1:b}}(\bm{x}^{b\leftarrow a}_{b-1})^\mathsf{T} \nonumber \\
    & \hspace{20pt} \circ \left[\bm{W}_{N}^\mathsf{T} \circ \bm{\sigma}_N' \circ (\bm{W}_{N} \cdot + \bm{b}_{N}) \right] \nonumber\\
    &\hspace{30pt} \circ T_{N-1:b} \circ (\Id+ T_{b-1:a+1}) \circ T_{a:1} (\bm{x}_0^{b\leftarrow a}) \notag\\
    & = \Big[ \big(\mathrm{J}_{T_{N-1:b}} ( \bm{x}^{b\leftarrow a}_{ b-1})\mathrm{J}_{T_{a:1}}(\bm{x}^{b\leftarrow a}_{ 0})\big)^\mathsf{T}\notag \\
    & \hspace{20pt}+\big(\mathrm{J}_{T_{N-1:b}}(\bm{x}^{b\leftarrow a}_{b-1}) \mathrm{J}_{T_{b-1:a+1}}(\bm{x}^{b\leftarrow a}_{a}) \mathrm{J}_{T_{a:1}}(\bm{x}^{b\leftarrow a}_{0}) \big)^\mathsf{T} \Big] \nonumber\\
    & \hspace{20pt} \circ \left[\bm{W}_{N}^\mathsf{T} \circ \bm{\sigma}_N' \circ (\bm{W}_{N} \cdot + \bm{b}_{N}) \right] \nonumber\\
    & \hspace{30pt} \circ T_{N-1:b} \circ (\Id + T_{b-1:a+1}) \circ T_{a:1} (\bm{x}_0^{b\leftarrow a}) \notag\\
    %
    %
    &= \left[ \prod_{\substack{n=1\\n\notin\{a+1, a+2, \cdots, b-1\}}}^{N-1}\mathrm{J}_{T_n}(\bm{x}_n^{b\leftarrow a})^\mathsf{T} \right] \circ \bm{W}_N^\mathsf{T}\tsup{T}_N(\bm{x}_0^{b\leftarrow a})
    \notag \\
    &\hspace{20pt}+ \left[ \prod_{n=1}^{N-1}\mathrm{J}_{T_n}(\bm{x}_n^{b\leftarrow a})^\mathsf{T} \right] \circ \bm{W}_N^\mathsf{T} \tsup{T}_N(\bm{x}_0^{b\leftarrow a})
    \notag \\
    &=\widetilde{R}^{b\leftarrow a}_{ 1}(\bm{x}_0^{b\leftarrow a}) +  \widetilde{R}_{1}(\bm{x}_0^{b\leftarrow a}).
\end{align}
The denoiser $D^{b\leftarrow a}$ is Lipschitz continuous because it can be written as a sum of $\widetilde{R}^{b\leftarrow a}_1$ and $\widetilde{R}_1$ which are Lipschitz continuous, as shown below.
The Lipschitz continuity of $\widetilde{R}^{b\leftarrow a}_1$ and $\widetilde{R}_1$ can actually be shown in an analogous way to Theorem~\ref{theorem_mol_NN}(b).
Indeed, the boundedness of $\tsup{t}_{n, i}$ follows from assumption \ref{assumption_lip_sigma}, and the Lipschitz continuity of $\tsup{T}_n$ is ensured by assumptions \ref{assumption_lip_sigma} and \ref{assumption_lip_sigmaprime} together with the Lipschitz continuity of $\Id+T_{b-1:a+1}$.
Thus, letting $F:=\widetilde{R}_n$, $R:=\widetilde{R}_{n+1}$, $T:=\tsup{T}_n$, and $\bm{W}:=\bm{W}_{n}^\mathsf{T}$ in Lemma~\ref{memo_lip2} verifies Lipschitz continuity of $\widetilde{R}_n$, as in Theorem~\ref{theorem_mol_NN}(b).
Likewise, the Lipschitz continuity of $\widetilde{R}^{b\leftarrow a}_n$ can be verified by letting $F:=\widetilde{R}^{b\leftarrow a}_n$, $R:=\widetilde{R}_{n+1}^{b\leftarrow a}$, $T:=\tsup{T}_n$, and $\bm{W}:=\bm{W}_{n}^\mathsf{T}$ for $n\in \{ 1,2,\cdots, a, b, \cdots, N \}$.
This completes the proof of the Lipschitz continuity of $D^{b\leftarrow a}$.
Since Theorem\ref{theorem_skip}(a) verifies the convexity of $\sum_{i=1}^{d_N} t_i^{b \leftarrow a}$, $D^{b \leftarrow a} = \nabla \left( \sum_{i=1}^{d_N} t_i^{b \leftarrow a} \right)$ is a MoL-Grad denoiser.
\qed
\ifCLASSOPTIONcaptionsoff
  \newpage
\fi

\bibliographystyle{IEEEtran}
\bibliography{refs}

\begin{thebibliography}{10}
\providecommand{\url}[1]{#1}
\csname url@samestyle\endcsname
\providecommand{\newblock}{\relax}
\providecommand{\bibinfo}[2]{#2}
\providecommand{\BIBentrySTDinterwordspacing}{\spaceskip=0pt\relax}
\providecommand{\BIBentryALTinterwordstretchfactor}{4}
\providecommand{\BIBentryALTinterwordspacing}{\spaceskip=\fontdimen2\font plus
\BIBentryALTinterwordstretchfactor\fontdimen3\font minus \fontdimen4\font\relax}
\providecommand{\BIBforeignlanguage}[2]{{%
\expandafter\ifx\csname l@#1\endcsname\relax
\typeout{** WARNING: IEEEtran.bst: No hyphenation pattern has been}%
\typeout{** loaded for the language `#1'. Using the pattern for}%
\typeout{** the default language instead.}%
\else
\language=\csname l@#1\endcsname
\fi
#2}}
\providecommand{\BIBdecl}{\relax}
\BIBdecl

\bibitem{yukawa2025_molgrad}
M.~Yukawa and I.~Yamada, ``Monotone {L}ipschitz-gradient denoiser: {E}xplainability of operator regularization approaches free from {L}ipschitz constant control,'' \emph{IEEE Trans.~Signal Process.}, vol.~73, pp. 3378--3393, 2025.

\bibitem{RUDIN1992_TV}
L.~I. Rudin, S.~Osher, and E.~Fatemi, ``Nonlinear total variation based noise removal algorithms,'' \emph{Physica D: Nonlinear Phenomena}, vol.~60, no.~1, pp. 259--268, 1992.

\bibitem{KHORAMIAN2012_penaltynorm}
S.~Khoramian, ``An iterative thresholding algorithm for linear inverse problems with multi-constraints and its applications,'' \emph{Applied and Computational Harmonic Analysis}, vol.~32, no.~1, pp. 109--130, 2012.

\bibitem{jian_nonconvexTV}
J.~Zou, M.~Shen, Y.~Zhang, H.~Li, G.~Liu, and S.~Ding, ``Total variation denoising with non-convex regularizers,'' \emph{IEEE Access}, vol.~7, pp. 4422--4431, 2019.

\bibitem{Rowl}
T.~Sasaki, Y.~Bandoh, and M.~Kitahara, ``Sparse regularization based on reverse ordered weighted ${L}_1$-norm and its application to edge-preserving smoothing,'' in \emph{Proc.~IEEE ICASSP}, 2024, pp. 9531--9535.

\bibitem{Ochs_nonconv_regu}
P.~Ochs, A.~Dosovitskiy, T.~Brox, and T.~Pock, ``On iteratively reweighted algorithms for nonsmooth nonconvex optimization in computer vision,'' \emph{SIAM J. Imag. Sci.}, vol.~8, no.~1, pp. 331--372, 2015.

\bibitem{Ben_nonconvexTV}
B.~Wang and Z.~Lan, ``Nonconvex high-order tv and $ \ell_0 $ norm-based method for image restoration,'' \emph{Electronic Research Archive}, vol.~33, pp. 3431--3449, 2025.

\bibitem{ven2013_PnP}
S.~V. Venkatakrishnan, C.~A. Bouman, and B.~Wohlberg, ``Plug-and-play priors for model based reconstruction,'' in \emph{Proc.~IEEE Global Conf. on Signal Inf. Process.}, 2013, pp. 945--948.

\bibitem{K_SVD_origin}
M.~Aharon, M.~Elad, and A.~Bruckstein, ``K-svd: An algorithm for designing overcomplete dictionaries for sparse representation,'' \emph{IEEE Trans. Signal Process.}, vol.~54, no.~11, pp. 4311--4322, 2006.

\bibitem{K_SVD}
M.~Elad and M.~Aharon, ``Image denoising via sparse and redundant representations over learned dictionaries,'' \emph{IEEE Trans. Image Process.}, vol.~15, no.~12, pp. 3736--3745, 2006.

\bibitem{BM3D}
K.~Dabov, A.~Foi, V.~Katkovnik, and K.~Egiazarian, ``Image denoising by sparse {3}-{D} transform-domain collaborative filtering,'' \emph{IEEE Trans. Image Process.}, vol.~16, no.~8, pp. 2080--2095, 2007.

\bibitem{TV_for_Pnp}
T.~Goldstein and S.~Osher, ``The split bregman method for {L}1-regularized problems,'' \emph{SIAM J. Imaging Sci.}, vol.~2, no.~2, pp. 323--343, 2009.

\bibitem{DnCNN}
K.~Zhang, W.~Zuo, Y.~Chen, D.~Meng, and L.~Zhang, ``Beyond a gaussian denoiser: Residual learning of deep {CNN} for image denoising,'' \emph{IEEE Trans. Image Process.}, vol.~26, no.~7, pp. 3142--3155, 2017.

\bibitem{IRCNN}
K.~Zhang, W.~Zuo, S.~Gu, and L.~Zhang, ``Learning deep {CNN} denoiser prior for image restoration,'' in \emph{IEEE CVPR}, 2017, pp. 2808--2817.

\bibitem{FFDNet}
K.~Zhang, W.~Zuo, and L.~Zhang, ``{FFDN}et: Toward a fast and flexible solution for {CNN}-based image denoising,'' \emph{IEEE Trans. Image Process.}, vol.~27, no.~9, pp. 4608--4622, 2018.

\bibitem{FFD_for_PnP}
X.~Yuan, Y.~Liu, J.~Suo, and Q.~Dai, ``Plug-and-play algorithms for large-scale snapshot compressive imaging,'' in \emph{Proc.~IEEE CVPR}, 2020, pp. 1444--1454.

\bibitem{Zang2022_drunet}
K.~Zhang, Y.~Li, W.~Zuo, L.~Zhang, L.~Van~Gool, and R.~Timofte, ``Plug-and-play image restoration with deep denoiser prior,'' \emph{IEEE Trans. Pattern Analysis and Machine Intelligence}, vol.~44, no.~10, pp. 6360--6376, 2022.

\bibitem{diffusion_PnP_IR}
Y.~Zhu, K.~Zhang, J.~Liang, J.~Cao, B.~Wen, R.~Timofte, and L.~V. Gool, ``Denoising diffusion models for plug-and-play image restoration,'' in \emph{Proc.~IEEE CVPRW}, 2023, pp. 1219--1229.

\bibitem{sreehari2016}
S.~Sreehari, S.~V. Venkatakrishnan, B.~Wohlberg, G.~T. Buzzard, L.~F. Drummy, J.~P. Simmons, and C.~A. Bouman, ``Plug-and-play priors for bright field electron tomography and sparse interpolation,'' \emph{IEEE Trans. Computat. Imag.}, vol.~2, no.~4, pp. 408--423, 2016.

\bibitem{romano2016_red}
Y.~Romano, M.~Elad, and P.~Milanfar, ``The little engine that could: Regularization by denoising ({RED}),'' \emph{SIAM J. Imaging Sci.}, vol.~10, no. 139, 11 2016.

\bibitem{sun2019}
Y.~Sun, B.~Wohlberg, and U.~S. Kamilov, ``An online plug-and-play algorithm for regularized image reconstruction,'' \emph{IEEE Trans. Computat. Imag.}, vol.~5, no.~3, pp. 395--408, 2019.

\bibitem{teodoro2019}
A.~M. Teodoro, J.~M. Bioucas-Dias, and M.~A.~T. Figueiredo, ``A convergent image fusion algorithm using scene-adapted gaussian-mixture-based denoising,'' \emph{IEEE Trans. Image Process.}, vol.~28, no.~1, pp. 451--463, 2019.

\bibitem{pesquet2021_MMO}
J.-C. Pesquet, A.~Repetti, M.~Terris, and Y.~Wiaux, ``Learning maximally monotone operators for image recovery,'' \emph{SIAM J. Imaging Sci.}, vol.~14, no.~3, pp. 1206--1237, 2021.

\bibitem{Nair2021}
P.~Nair, R.~G. Gavaskar, and K.~N. Chaudhury, ``Fixed-point and objective convergence of plug-and-play algorithms,'' \emph{IEEE Trans. Computat. Imag.}, vol.~7, p. 337–348, 2021.

\bibitem{sun2021}
Y.~Sun, Z.~Wu, X.~Xu, B.~Wohlberg, and U.~S. Kamilov, ``Scalable plug-and-play admm with convergence guarantees,'' \emph{IEEE Trans. Computat. Imag.}, vol.~7, pp. 849--863, 2021.

\bibitem{RED_pro}
R.~Cohen, M.~Elad, and P.~Milanfar, ``Regularization by denoising via fixed-point projection ({RED}-{PRO}),'' \emph{SIAM J. Imaging Sci.}, vol.~14, no.~3, pp. 1374--1406, 2021.

\bibitem{reehost2019}
E.~T. Reehorst and P.~Schniter, ``Regularization by denoising: Clarifications and new interpretations,'' \emph{IEEE Trans. Computat. Imag.}, vol.~5, no.~1, pp. 52--67, 2019.

\bibitem{gradient_step_origin}
R.~Cohen, Y.~Blau, D.~Freedman, and E.~Rivlin, ``It has potential: Gradient-driven denoisers for convergent solutions to inverse problems,'' in \emph{Advances in Neural Information Processing Systems}, vol.~34.\hskip 1em plus 0.5em minus 0.4em\relax Curran Associates, Inc., 2021, pp. 18\,152--18\,164.

\bibitem{hurault2022_gradient}
S.~Hurault, A.~Leclaire, and N.~Papadakis, ``Gradient step denoiser for convergent plug-and-play,'' in \emph{Proc.~ICLR}, 2022.

\bibitem{hurault2022proximal}
------, ``Proximal denoiser for convergent plug-and-play optimization with nonconvex regularization,'' in \emph{Proc. ICML}, ser. Proceedings of Machine Learning Research, K.~Chaudhuri, S.~Jegelka, L.~Song, C.~Szepesvari, G.~Niu, and S.~Sabato, Eds., vol. 162.\hskip 1em plus 0.5em minus 0.4em\relax PMLR, 17--23 Jul 2022, pp. 9483--9505.

\bibitem{ryu2019_nonexptrain_by-norm}
E.~Ryu, J.~Liu, S.~Wang, X.~Chen, Z.~Wang, and W.~Yin, ``Plug-and-play methods provably converge with properly trained denoisers,'' in \emph{Proc.~ISML}, ser. Proceedings of Machine Learning Research, vol.~97.\hskip 1em plus 0.5em minus 0.4em\relax PMLR, 2019, pp. 5546--5557.

\bibitem{tsuzuku2018}
Y.~Tsuzuku, I.~Sato, and M.~Sugiyama, ``Lipschitz-margin training: scalable certification of perturbation invariance for deep neural networks,'' ser. NIPS'18.\hskip 1em plus 0.5em minus 0.4em\relax Red Hook, NY, USA: Curran Associates Inc., 2018, p. 6542–6551.

\bibitem{pesquet_icassp20}
A.~Neacsu, J.-C. Pesquet, and C.~Burileanu, ``Accuracy-robustness trade-off for positively weighted neural networks,'' in \emph{Proc. ICASSP}, 2020, pp. 8389--8393.

\bibitem{Numerical_Optimization}
S.~J.~W. Jorge~Nocedal, \emph{Numerical Optimization}, ser. Springer Series in Operations Research and Financial Engineering.\hskip 1em plus 0.5em minus 0.4em\relax Springer New York, 2006.

\bibitem{nguyen_nonnegativeRBM}
T.~Nguyen, T.~Tran, D.~Phung, and S.~Venkatesh, ``Learning parts-based representations with nonnegative restricted boltzmann machine,'' \emph{J. Machine Learn. Research}, vol.~29, p. To Appear, 01 2013.

\bibitem{hosseini_nonnegative_aut}
E.~Hosseini-Asl, J.~M. Zurada, and O.~Nasraoui, ``Deep learning of part-based representation of data using sparse autoencoders with nonnegativity constraints,'' \emph{IEEE Trans. Neural Networks and Learning Systems}, vol.~27, no.~12, pp. 2486--2498, 2016.

\bibitem{convex_analysis}
H.~H. Bauschke and P.~L. Combettes, \emph{Convex Analysis and Monotone Operator Theory in Hilbert Spaces}, 2nd~ed.\hskip 1em plus 0.5em minus 0.4em\relax New York, NY: Springer, 2017.

\bibitem{donoho1995_softshrink}
D.~Donoho, ``De-noising by soft-thresholding,'' \emph{IEEE Trans. Information Theory}, vol.~41, no.~3, pp. 613--627, 1995.

\bibitem{Fan2001_SCAD_shrink}
J.~Fan and R.~Li, ``Variable selection via nonconcave penalized likelihood and its oracle properties,'' \emph{J. the American Statistical Association}, vol.~96, no. 456, pp. 1348--1360, 2001.

\bibitem{chartrand2007_lp}
R.~Chartrand, ``Exact reconstruction of sparse signals via nonconvex minimization,'' \emph{IEEE Signal Process. Lett.}, vol.~14, no.~10, pp. 707--710, 2007.

\bibitem{candes2007_reweighted_l1}
E.~J. Candes, M.~B. Wakin, and S.~P. Boyd, ``Enhancing sparsity by reweighted $\ell_1$ minimization,'' \emph{J. Fourier Analysis and Applications}, vol.~14, pp. 877--905, 10 2008.

\bibitem{Zhang2010_MCP_shrink}
C.-H. Zhang, ``Nearly unbiased variable selection under minimax concave penalty,'' \emph{The Annals of Statistics}, vol.~38, 02 2010.

\bibitem{yao2017_family_nonconvex}
Q.~Yao and J.~T. Kwok, ``Efficient learning with a family of nonconvex regularizers by redistributing nonconvexity,'' \emph{J. Machine Learn. Research}, vol.~18, no.~1, p. 6574–6625, Jan. 2017.

\bibitem{yukawa2023_limes}
M.~Yukawa, H.~Kaneko, K.~Suzuki, and I.~Yamada, ``Linearly-involved moreau-enhanced-over-subspace model: Debiased sparse modeling and stable outlier-robust regression,'' \emph{IEEE Trans. Signal Process.}, vol.~71, pp. 1232--1247, 2023.

\bibitem{suzuki_external}
K.~Suzuki and M.~Yukawa, ``External division of two proximity operators: An application to signal recovery with structured sparsity,'' in \emph{Proc.~IEEE ICASSP}, 2024, pp. 9471--9475.

\bibitem{combettes_lip}
P.~L. Combettes and J.-C. Pesquet, ``{L}ipschitz certificates for layered network structures driven by averaged activation operators,'' \emph{SIAM Journal on Mathematics of Data Science}, vol.~2, no.~2, pp. 529--557, 2020.

\bibitem{Xiao2017_fMNIST}
H.~Xiao, K.~Rasul, and R.~Vollgraf, ``Fashion-mnist: a novel image dataset for benchmarking machine learning algorithms,'' \emph{ArXiv}, vol. abs/1708.07747, 2017.

\bibitem{BSDS}
D.~Martin, C.~Fowlkes, D.~Tal, and J.~Malik, ``A database of human segmented natural images and its application to evaluating segmentation algorithms and measuring ecological statistics,'' in \emph{Proc.~8th Int'l Conf. Computer Vision}, vol.~2, July 2001, pp. 416--423.

\bibitem{Levin_blurkernel}
A.~Levin, Y.~Weiss, F.~Durand, and W.~T. Freeman, ``Understanding and evaluating blind deconvolution algorithms,'' in \emph{Proc.~IEEE CVPR}, 2009, pp. 1964--1971.

\bibitem{Flicker2K}
B.~Lim, S.~Son, H.~Kim, S.~Nah, and K.~M. Lee, ``Enhanced deep residual networks for single image super-resolution,'' in \emph{Proc.~IEEE CVPRW}, 2017, pp. 1132--1140.

\bibitem{Agustsson2017_DIV2K}
E.~Agustsson and R.~Timofte, ``Ntire 2017 challenge on single image super-resolution: Dataset and study,'' in \emph{Proc.~IEEE CVPRW}, July 2017.

\bibitem{ma2017_waterloo}
K.~Ma, Z.~Duanmu, Q.~Wu, Z.~Wang, H.~Yong, H.~Li, and L.~Zhang, ``{Waterloo Exploration Database}: New challenges for image quality assessment models,'' \emph{IEEE Trans. Image Process.}, vol.~26, no.~2, pp. 1004--1016, Feb. 2017.

\bibitem{Kingma2014_Adam}
D.~P. Kingma and J.~Ba, ``Adam: A method for stochastic optimization,'' \emph{CoRR}, vol. abs/1412.6980, 2014.

\end{thebibliography}

%


\end{document}